\documentclass[sn-mathphys,Numbered]{sn-jnl}
\jyear{2023}
\usepackage{graphicx}%
\usepackage{multirow}%
\usepackage{amsmath,amssymb,amsfonts}%
\usepackage{amsthm}%
\usepackage{mathrsfs}%
\usepackage[title]{appendix}%
\usepackage{caption}
\theoremstyle{thmstyleone}
\newtheorem{theorem}{Theorem}[section]
\newtheorem{proposition}[theorem]{Proposition}
\theoremstyle{thmstyletwo}%
\newtheorem{example}{Example}%
\newtheorem{remark}{Remark}%
\newtheorem{lemma}{Lemma}%
\newtheorem{corollary}{Corollary}%
\theoremstyle{thmstylethree}%
\begin{document}
\title{New copula families and mixing properties}          
\author*[1]{\fnm{Martial} \sur{Longla}}\email{mlongla@olemiss.edu}
\affil*[1]{\orgdiv{Department of mathematics}, \orgname{University of Mississippi}, \orgaddress{\street{University Ave}, \city{University}, \postcode{38677}, \state{Mississippi}, \country{US}}}

\abstract{We characterize absolutely continuous symmetric copulas with square integrable densities in this paper. This characterization is used to create new copula families, that are perturbations of the independence copula. The full study of mixing properties of Markov chains generated by these copula families is conducted. An extension that includes the Farlie-Gumbel-Morgenstern family of copulas is proposed. We propose some examples of copulas that generate non-mixing Markov chains, but whose convex combinations generate $\psi$-mixing Markov chains. Some general results on $\psi$-mixing are given. The Spearman's correlation $\rho_S$ and Kendall's $\tau$ are provided for the created copula families.  Some general remarks are provided for $\rho_S$ and $\tau$. A central limit theorem is provided for parameter estimators in one example. A simulation study is conducted to support derived asymptotic distributions for some examples.}


\maketitle

\section{Introduction}
A bivariate copula, which is defined as a restriction to the unit square $[0; 1]^2$ of a bivariate joint cumulative distribution function with uniform marginals, is a tool that has been gaining popularity in dependence modeling. An equivalent definition can be found in Nelsen (2006) \cite{N}, where one can find several related notions. The popularity of copulas is due to Sklar's theorem that relates them to the joint distribution and marginals of a bivariate random vector (see Sklar (1959) \cite{SK}).

Many authors have worked on building copulas with various properties. Some constructions can be found in Nelsen (2006) \cite{N}, where several construction methods are presented. Arakelian and Karlis (2014)\cite{AK} studied mixtures of copulas, that were investigated for mixing by Longla (2015)\cite{Lo2015}.   Longla (2014) \cite{L2} constructed copulas based on prescribed $\rho$-mixing properties obtained from extending results of Beare (2010)\cite{Beare} on mixing for copula-based Markov chains. Longla et al. (2022) \cite{L0} and Longla et al. (2022b)\cite{L00} constructed copulas in general based on perturbations. In these two works, two perturbation methods were considered. For the first method, copulas are perturbed at the level of the variables by adding some noise components to each of the variables.  For the second method, a copula  $C(u,v)$ is perturbed by creating $C_D(u,v)=C(u,v)+D(u,v)$ and requiring that $C_D(u,v)$ satisfies the definition of a copula. This last method of perturbation was also presented in Komornik et al (2017) \cite{KK} and the references therein. Several other authors have considered extensions of existing copula families via other methods, that are not the focus of this work, see for instance Morillas (2005)\cite{Morillas}, Aas et al (2009) \cite{AAS}, Klement et al (2017) \cite{KL}, Chesneau (2021) \cite{Chesneau} among others. 

Chesneau (2021)\cite{Chesneau} considered multivariate trigonometric copulas, which seem close to one of the examples of copula families of this paper, but are in fact very far from the work we provide here. We are concerned with the important question of characterizing absolutely continuous symmetric copulas with square integrable densities. After characterizing such copulas, we extract some general copula families with specific functions, one of which is made of trigonometric functions. These copula families come with pre-difined mixing structure of the copula-based Markov chains they generate. This, in general, can not be said about the copula families of Chesneau (2021)\cite{Chesneau}, which where published when this paper was already written. We also provide the mixing structure of Markov chains generated by any of the copulas constructed in this paper.  Therefore, one of the main points of the paper is the presentation of several copula families and the mixing structure of the Markov chains they generate. Moreover, for each of these copulas, we provide joint distributions of any two variables along the Markov chain they generate; and show that the copulas of these variables also belong to the initial copula family. We show that the set of all copulas for each of the selected basis of functions is closed under the operation $*$, defined by
$C^1(u,v)=C(u,v),$ for $n>1$ \[ C^n(u,v)=C^{n-1}*C(u,v)=\int_0^1 C^{n-1}_{,2}(u,t)C_{,1}(t,v)dt;\]
where $C_{,i}(u,v)$ is the derivative with respect to the $i^{th}$ component of $(u,v)$.  This product plays a central role in our study and its properties can be found in Darsow et al. (1992) \cite{DAR}. Recall that for $n=2$ this product is the joint distribution of $(X_1,X_3)$ when $(X_1,X_2,X_3)$ is a stationary Markov chain with copula $C(u,v)$ and the uniform marginal distribution. In general, the copula $C^{n}(u,v)$ is the joint distribution of $(X_0, X_n)$ when $X_0, \cdots, X_n$ is a stationary Markov chain generated by $C(u,v)$ and the Uniform(0,1) distribution. This is a simple consequence of applying recursively the product formula (see Longla (2022) \cite{L0}). 

The product $*$ and related notions are used in this paper to derive properties of Markov chains, including mixing and association. In Longla et al. (2022) \cite{L0}, measures of association such as Spearman's $\rho$ and Kendall's $\tau$ were studied for perturbations of copulas. Recall that in terms of copulas, these coefficients are defined as 
\begin{equation} \rho_S(C)=12\int_0^1\int_0^1C(u,v)dudv-3,\label{Rho}\end{equation} \begin{equation} \text{and}\quad \tau(C)=1-4\int_0^1\int_0^1C_{,1}(u,v)C_{,2}(u,v)dudv. \label{Ken}
\end{equation}
 Formulas \eqref{Rho} and \eqref{Ken} are respectively formula (5.1.15c) and formula (5.1.12) of Nelsen (2006)\cite{N}. For definitions and practical use of these measures, we refer the reader to Nelsen (2006) \cite{N} and references therein. In this paper, our interest for these measures is in finding the relationship between them and the parameters of the copulas, then propose an estimation procedure. 
Among copula families, for which we derive measures of association, is the family of copulas based on Legendre polynomials. This family is an extension of the Farlie-Gumbel-Morgenstern (FGM) copula family. The FGM has been extensively studied in the literature with several extensions (see Hurlimann (2017)\cite{H} or Ebaid et al. (2022) \cite{Eb}). Our extension is unrelated to extensions that exist in the literature. We use a new approach to extend this family in a way that has not been done in the literature. Moreover, we show that copulas from this family generate $\psi$-mixing stationary Markov chains for all values of the family parameter $\theta\in[-1,1]$. This improves a previous result of Longla et al (2022) \cite{L0}, who failed to provide the proof for the boundary points $\theta=\pm 1$. The new characterization of copulas obtained in this paper resulted in a new elegant proof of $\psi$-mixing that includes the boundary points.
 For more on mixing coefficients, see Bradley (2005)\cite{Bradley}, Beare (2010) \cite{Beare} and the references therein. To avoid confusions in this paper, we say that a function $f(x)$ is a version of the function $g(x)$ if and only if the Lebesgue measure (in the appropriate dimension) of the set on which they differ is 0.

This paper is structured as follows. In section 2, we present a characterization of copulas with square integrable densities. The obtained characterization is used to construct new copula families. Among these general functional families of copulas, are finite and infinite sums, including new sine and cosine copula families. A second group of copulas that is shown as example here is based on Legendre polynomials and contains the FGM copula. We also provide the Spearman's $\rho_S$ and Kendall's $\tau$ coefficients for each of the examples of copula families along with some useful remarks. All these results are new to the best of our knowledge. In section 3, we develop mixing properties of the extended FGM copula family based on Legendre polynomials. Several examples are provided here along with some copulas that don't generate mixing Markov chains, but have convex combinations that generate Mixing Markov chains. These examples answer an open question on convex combinations of copulas and mixing. A new general result on mixing for Markov chains generated by copulas with square integrable densities is given. 
We provide an extension of a previous result of Longla et al (2022c) \cite{LoMous} on $\psi$-mixing and $\psi'$-mixing. Namely, 
we show that $\psi$-mixing follows from the two conditions $c^{s_1}(u,v)>0$ on a set of Lebesgue measure 1 and $c^{s_2}(u,v)<2$ for all $(u,v)\in [0,1]^2$. In section 4 we provide some parameter estimators and their asymptotic distributions via centra limit theorems. We provide here a simulation study. Section 5 covers comments and conclusions. The Appendix of proofs ends the paper.

\section{Square integrable symmetric densities}
In this section, we provide a characterization of copulas with square integrable densities. The aim is to show that any such copula can be represented as a sum of possibly non-zero infinitely many terms. Such a representation indicates a new way to construct new copula families by modifying the independence copula $\Pi(u,v)=uv$. We use this method in this section to provide new sets of copulas and compute their coefficients $\rho_S$ and $\tau$. 
 
Recall that a copula in general is defined as a multivariate cumulative distribution function, for which the restriction on $[0,1]^2$ has uniform marginals on (0,1). For any bivariate random variable $(X,Y)$ with continuous marginal distributions $(F_X, F_Y)$, the copula of $(X,Y)$ is the joint cumulative distribution function of $(F_X(X), F_Y(Y))$. Note that both $F_X(X)$ and $F_Y(Y)$ are uniform on $[0,1]$. This is why most works on copula theory ignore the marginal distributions and are limited to uniform marginals. This is done without loss of generality, when variables are assumed continuous (see Nelsen (2006) \cite{N}, Durante and Sempi (2016) \cite{DS} or Sklar (1959) \cite{SK} for more on the topic). This work is concerned mainly with absolutely continuous symmetric bivariate copulas $C(u,v)$ for which the density function $c(u,v)$ exist and satisfies 
\begin{equation}\label{AC}
C(u,v)=\int_0^u\int_0^vc(s,t)dtds.
\end{equation}
Formula \eqref{AC} follows from the fact that the singular part of the copula vanishes for absolutely continuous copulas (see Darsow et al. (1992) \cite{DAR}). It is known that the density of a copula is a positive bivariate function on $[0,1]^2$. Many authors have worked on copulas and their properties. Copulas are used to model dependence in various applied fields for  problems that include but are not limited to estimation, classification and statistical tests of hypotheses. When used to model dependence, they help establish mixing properties of Markov chains, that are useful in establishing central limit theorems for sample averages of functions of observations (see Doukhan et al. (2009)\cite{Doukhan} Chen and Fan (2006)\cite{Chen}, Peligrad and Utev (1997)\cite{PelUt}, and Jones (2004) \cite{Jones}). For more on properties of copulas, see Darsow et al. (1992) \cite{DAR} and Beare (2010) \cite{Beare}. 

In this paper, we say that a copula $C(u,v)$ is symmetric if and only if $C(u,v)=C(v,u)$ for all $(u,v)\in[0,1]^2$. Such copulas are useful in modeling reversible Markov chains. For a symmetric absolutely continuous copula, the density $c(u,v)$ is square integrable if and only if 
\begin{equation} \int_0^1\int_0^1 [c(s,t)]^2dtds<\infty.
\end{equation}
Square integrable copula densities are kernels of Hilbert-Shmidt operators on the Hilbert space $L^2(0,1)$. This means, that the linear operator defined by 
\begin{equation}
Kf(x)=\int_0^1c(x,y)f(y)dy \label{K}
\end{equation}
is a Hilbert-Shmidt operator in $L^2(0,1)$ (see Ferreira and Manegatto (2009) \cite{Fer}). Let $\{\varphi_k(x), k\in\mathbb{N}\}$ be an orthonormal basis made of eigenfunctions of the operator $K$. $K\varphi_k(x)=\lambda_k\varphi_k(x),$ where $\lambda_k$ is called eigenvalue of $K$ associated to the eigenfunction $\varphi_k(x)$ and the equality is given in the sense of $L^2(0,1)$. Based on this, for a symmetric copula with square integrable density, we have (see also Beare (2010) \cite{Beare}, Longla (2014) \cite{L2} ):
\begin{itemize}
\item The sequence of $\lvert{}\lambda_k\lvert{}$ has a finite number of values including their multiplicities, or converges to $0$.  
\item  $\{\varphi_k(x), k\in\mathbb{N}\}$ can be selected to form an orthonormal system in $L^2(0,1)$.
\end{itemize} 
Due to the fact that $L^2(0,1)$ is a separable Hilbert space and the orthonormal system made of eigenfunctions of this operator is complete, it follows that formula \eqref{cond1} holds in the sense of $L^2(0,1)$. Therefore, by Mercer's theorem (see formula 27 page 445 of Mercer (1909) \cite{Mercer} or Theorem 2.4 of Ferreira and Manegatto (2009) \cite{Fer})
\begin{equation} c(u,v)=\sum_{k=1}^{\infty}\lambda_k\varphi_k(u)\varphi_k(v), \label{cond1}
\end{equation}
where convergence of the series is point-wise and uniform on $(0,1)^2$ when the opertor is non-negative definite (when all eigenvalues are non-negative).  The same conclusion can be achieved when all eigenvalues are non-positive (or the operator is non-positive). In fact, Mercer (1909) \cite{Mercer} showed that for any continuous non-negative definite kernel $k(x,y)$, the series converges point-wise and uniformly. A non-negative definite kernel is a function $k(x,y)$ such that 
\[\sum_{i,j=1}^n k(x_i,x_j)s_is_j \ge 0 \quad \text{for all}\quad x_i,x_j, s_i, s_j.\]
Mercer (1909) \cite{Mercer} pointed that the above definition is equivalent to \[ \int\int k(x(t), x(s))x(t)x(s)dsdt\ge 0.\]
 Moreover, it is obvious that for a continuous symmetric kernel $k(x,y)$, the trace of the operator is \[\sum_{j=1}^\infty \lambda_j=\int_0^1k(x,x)dx.\]
Therefore, the series of partial sums of eigenvalues is convergent when the density of the copula is symmetric, square integrable, positive definite and continuous on $[0,1]^2$. It is important to mention that these conditions impose that eigenvalues be all positive (or all negative), therefore ruling out some of the copulas. But the decomposition can be extended by continuity to copulas with eigenvalues of both signs.

\begin{proposition} \label{pr1} For any absolutely continuous copula $C(u,v)$ the following holds.
\begin{enumerate}
\item The kernel operator \eqref{K} has and eigenvalue $\lambda=1$ associated to the eigenfunction $\varphi(x)=1$;
\item There exists a basis of eigenfunctions of \eqref{K} for any square integrable symmetric copula density. Moreover, there exists a decomposition \eqref{cond1} containing $\varphi_1(x)=1$ associated to the eigenvalue $\lambda_1=1$. 
\item A selected basis might not contain $\varphi(x)=1$ in the case when the eigenvalue $\lambda=1$ has multiplicity at least equal to $2$, generating a subspace of dimension higher than $1$ with orthonormal basis not containing $\varphi(x)=1$.
\end{enumerate}
\end{proposition}

Note that when the eigenvalue $\lambda=1$ has multiplicity $s>1$ in the decomposition, $s$ has to be finite because of square integrability.
It is also worth mentioning that it is impossible for some functions to be eigenfunctions associated to the eigen value $1$. For instance, $\varphi_k(x)$ has to be bounded.

We devote the following subsections of this work to the case when $\varphi(x)=1$ is used. This case by itself is interesting because it turns into the question of perturbations of the independence copula (see Longla et al. (2022) \cite{L0}). As we will see below, the mixing structure of Markov chains generated by our constructed copula families is easily established in general. This eases the study of estimators of coefficients and measures of association. 

\subsection{Perturbations of the copula $\Pi(u,v)$}
The copula $\Pi(u,v)=uv$ is called independence copula. It is used to model data under the assumption of independence. This means that the copula of random variables that are independent is $\Pi(u,v)$. The notion of perturbation has been recently used in many papers to address modifications that are done to a copula to improve estimation results. This is done by introducing some level of dependence in the case of the independence copula (see Durante et al (2013) \cite{F.S}) through a function $D(u,v)$. When $\Pi(u,v)$ is perturbed, the obtained copula exhibits some dependence that can be estimated to fit better the data at hands. 

Various types of perturbations have been considered in Durante et al (2013) \cite{F.S}, Komornik et al (2017) \cite{KK}, Longla et al (2022) \cite{L0}, Longla et al (2022c)\cite{LoMous} and the references therein. In this section, we investigate perturbations of the kind $C_D(u,v)=\Pi(u,v)+D(u,v)$, where $D(u,v)$ is the perturbation factor to be determined using \eqref{cond1}. This perturbation method was studied in Komornik et al (2017) \cite{KK} and the references therein. Longla et al (2022) \cite{L0} and Longla et al (2022c) \cite{LoMous} considered mixing for Markov chains generated by perturbations with $D(u,v)=\theta (\Pi(u,v)-C(u,v))$ and $D(u,v)=\theta (M(u,v)-C(u,v))$, where $M(u,v)=\min(u,v)$ is the  Fr\'{e}chet-Hoeffding upper bound. 

\begin{remark}
The use of $M(u,v)$ to perturb a copula $C(u,v)$ for data with uniform marginals assumes that values of the Markov chain have a positive probability of repeating themselves, before providing new values from the transition copula $C(u,v)$. 
\end{remark}

Considering the decomposition \eqref{cond1}, when the first term uses $\varphi(x)=1$, reindexing the remaining functions $\varphi_k(x)$ and integrating leads to 
\begin{equation}
C(u,v)=uv+\sum_{k=1}^{\infty} \lambda_k\int_{0}^u \varphi_k(x)dx\int_0^v\varphi_k(y)dy. \label{cop}
\end{equation}
This implies $C(u,v)=uv+D(u,v)$, where \[ D(u,v)=\sum_{k=1}^{\infty} \lambda_k\int_{0}^u \varphi_k(x)dx\int_0^v\varphi_k(y)dy.\]

Note that $D(u,0)=D(u,1)=0$ and $D(u,v)$ is symmetric. Therefore, the perturbation is a copula when its second order mixed derivative is non-negative (refer to Longla et al. (2022) \cite{L0}). So, any representation \eqref{cop} with $\lambda_k\in\mathbb{R}$ defines a copula with square integrable density if and only if
\begin{equation}
1+\sum_{k=1}^{\infty}\lambda_k\varphi_k(x)\varphi_k(y)\ge 0 \quad \mbox{for all $(x,y)\in [0,1]^2$ and} \label{cond}
\end{equation}
\begin{equation}
\sum_{k=1}^{\infty}\lambda_k^2 <\infty. \label{lambda}
\end{equation}
Condition \eqref{lambda} is required for square integrability of the density and takes into account the fact that $\varphi_k(x)$ form an orthonormal system of functions. Condition \eqref{cond} is satisfied when the following holds:
\begin{equation}
1+\sum_{k=1}^{\infty} \lambda_k \alpha_k \ge 0, \quad \alpha_k=\begin{cases} 
\max\varphi^2_k & \text{if } \quad \lambda_k < 0 \\
\min\varphi_k\max\varphi_k &\text{if} \quad \lambda_k > 0 .
\end{cases} \label{cond2}
\end{equation}
Note that condition \eqref{cond2} is not necessary. For different values of $k$, the functions $\varphi_k(x)$ might have their minima or maxima at different points. It is a strong condition, but easy to handle in practice. Note that when condition \eqref{cond2} holds, if $\varphi_k(x)$ are uniformly bounded, then \eqref{lambda} holds. This can be shown by first showing that the sum of $\lvert\lambda_k\lvert $ is convergent and larger than the sum in \eqref{lambda}. For all systems of functions that are considered in this section, all these conditions are satisfied. Thus, the following characterization theorems hold.

\begin{theorem}\label{Th1}
Any symmetric absolutely continuous copula with square integrable density can be represented by \eqref{cop}, where $\{1, \varphi_k(x), k>1, k\in\mathbb{N}\}$ is an orthonormal basis associated to the eigenvalues $\{1, \lambda_k, k>1, k\in\mathbb{N}\}$ of the operator \eqref{K}.
\end{theorem}
 In Theorem \ref{Th1}, equality is understood as equality of operators in $L^2(0,1)$. This equality turns into point-wise or uniform convergence when the density of the copula is continuous and defines a positive semi-definite operator. 
 
\begin{theorem}\label{T1} Assume that $(\lambda_k\in\mathbb{R}, k>1)$ is such that conditions \eqref{cond2} and \eqref{lambda} hold
for an orthonormal system of functions of $L^2(0,1)$ that contains $\varphi_1(x)=1$. The function \eqref{cop} is an absolutely continuous symmetric copula with square integrable density that we call Type-I Longla copula with notation $L(u,v, \Lambda, \Phi)$. 
\end{theorem}
For any system of functions $\varphi_k(x)$ that is uniformly bounded, condition \eqref{cond2} implies condition \eqref{lambda}.  Theorem \ref{T1} opens the ground for various copula families with functions as parameters. These copulas can have finite or infinite sums of terms. They depend both on the system of functions $\varphi_k(x)$ and the coefficients $\lambda_k$. The characterization formulated in the following general result summarizes the above analysis.

\subsection{Examples of new copula families}
 Theorem \ref{T1} and Theorem \ref{Th1} provide a relationship between extreme values of the functions $\varphi_k(x)$ and the possible eigenvalues $\lambda_k$ of \eqref{K}. For example, a function $\varphi(x)$ with maximum $2$ and minimum $-3$ cannot be associated to an eigenvalue less than $-1/9$ or larger than $1/6$, when all other coefficients are equal to 0. The larger the values of $\lvert{}\varphi_k(x)\lvert{}$, the smaller the maximal range of $\lambda_k$, when the sum \eqref{cond1} has only $\varphi_k(x)$ and $1$. Setting $\Phi(x)=\int_0^x\varphi(s)ds$ for any mean-zero normal function $\varphi(x)$ on $[0,1]$, we have the copula
\begin{equation} C(u,v)=uv+\lambda \Phi(u)\Phi(v), \quad \frac{-1}{\max \varphi^2(x)} \le \lambda \le \frac{-1}{\max\varphi(x)\min\varphi(x)}.\label{cop2} \end{equation}
 An example of copula \eqref{cop2}  is defined for
$ -\min(\frac{1}{\alpha},\alpha) \le \lambda \le1,$ $\alpha\in (0,\infty)$ by
\begin{equation}
C_{\alpha,\lambda}(u,v)= \begin{cases} (1+\lambda\alpha)uv \quad \text{if $0\le u,v\le\frac{1}{\alpha+1}$} \\ uv+\lambda u(1-v) \quad \text{if $0\le u\le \frac{1}{\alpha+1}<v\le 1$} \\ uv+\lambda v(1-u) \quad \text{if $0\le v\le \frac{1}{\alpha+1}<u\le 1$} \\ uv+\frac{\lambda}{\alpha}(1-u)(1-v), \quad \text{if $\frac{1}{\alpha+1}<u,v\le 1$}. \end{cases} \label{cops2}
\end{equation}
Copula \eqref{cops2} is also a representative of the only class of copulas with square integrable symmetric densities based on a single function $\varphi(x)$ that takes $2$ values on $[0,1]$. This class of copulas is defined using $$\varphi(x)=\sqrt{\alpha}\mathbb{I}(x\in A)-\frac{1}{\sqrt{\alpha}}\mathbb{I}(x\in [0,1]\setminus A),$$ where $A$ is a set of measure $\frac{1}{\alpha+1}$ for some strictly postive real number $\alpha$. Copula \eqref{cops2} is obtained via integration to get $\Phi(x)$ of formula \eqref{cop2} for  $$\varphi(x)=\sqrt{\alpha}\mathbb{I}(0\le x<\frac{1}{\alpha+1})-\frac{1}{\sqrt{\alpha}}\mathbb{I}(1\ge x\ge \frac{1}{\alpha+1}).$$ 
For copula \eqref{cops2}, while the parameter $\lambda$ answers for dependence (we will show later that its absolute value is the maximal coefficient of correlation), the parameter $\alpha$ answers for the portions of $[0,1]^2$ that have constant likelihood values. In fact, the density of this copula is constant on each of the four portions of its support. Via simple computations we establish that for  copula \eqref{cops2},  \[\rho_S(C_{\alpha,\lambda})=\frac{3\alpha\lambda}{(1+\alpha)^2} \quad \text{and}\quad \tau(C_{\alpha,\lambda})=\frac{2\alpha\lambda}{(1+\alpha)^2}.\] 

\begin{remark}
The coefficients $\rho_S(C_{\alpha,\lambda})$ and $\tau(C_{\alpha,\lambda})$ show that for the purpose of correlation or association, the range of $\alpha$ can be limited to $(0,1]$, because $1/\alpha$ and $\alpha$ produce the same value for $\rho_S(C_{\alpha,\lambda})$ and $\tau(C_{\alpha,\lambda})$. The maximum of these coefficients is reached at $\alpha=1$. 
\end{remark}
Interest in the copula family $C_{1,\lambda}(u,v)$ is in the fact that the range of the Pearson correlation coefficient of functions of variables $(U,V)$ generated by this copula family is $[-1,1]$. For $\lambda\ne\pm1$, the procedure to generate a stationary Markov chain from $C_{1,\lambda}(u,v)$ with the uniform marginal distribution is as follows. We generate an observation $U_0$ from $Uniform(0,1)$; then for every integer $i$, generate $Q_i$ from $Uniform(0,1)$ and define $U_i$ using the quantile function of the conditional distribution of $U_i\lvert{}U_{i-1}$. The quantile function is the inverse with respect to $v$ of the derivative with respect to $u$ of $C(u,v)$ at the point $(U_{i-1},U_i)$ (see Longla and Peligrad (2012) \cite{LoMag}). We obtain the formula
\begin{equation}\label{wl}
U_i=\begin{cases} \frac{Q_i}{1+\lambda}, \quad \text{if}\quad 0<U_{i-1}<.5,\quad Q_i < .5+.5\lambda \\
\frac{Q_i-\lambda}{1-\lambda}, \quad \text{if}\quad 0<U_{i-1}<.5,\quad Q_i > .5+.5\lambda\\
\frac{Q_i}{1-\lambda}, \quad \text{if}\quad 1>U_{i-1}>.5,\quad Q_i < .5-.5\lambda\\
\frac{Q_i+\lambda}{1+\lambda}, \quad \text{if}\quad 1>U_{i-1}>.5,\quad Q_i > .5-.5\lambda.
\end{cases}
\end{equation}
Recall that the joint cumulative distribution function of $(U_0, U_n)$ is the copula $C_{1, \lambda^n}(u,v)$. Therefore, the following holds.
\begin{lemma} For any stationary Markov chain generated by $C_{1,\lambda}(u,v)$ and the uniform marginal distribution for $\lambda\ne \pm1$ , 
we have $U_n-Q_n\to 0$ in probability when $n\to \infty$.
\end{lemma}
The limiting value of $U_n$ as $n\to\infty$ doesn't depend on the initial value $U_0$. This is true because $\lambda$ is replaced by $\lambda^n\to 0$ in formula \eqref{wl} and $Q_n$ is a random sample obtained independently. This fact is one of the selling points of the study of mixing properties. The intial point of the Markov chain doesn't affect the long run behavior under mixing assumptions.

An extension of the copula family \eqref{cops2} is obtained considering for all $s\in\mathbb{N^*}$, $\delta_i=\mathbb{I}(a_i\le u, v<a_{i+1})$ and $\lvert \theta_i\lvert\leq 1$. This extension has density
\begin{eqnarray}
c(u,v)=1+\sum_{i=1}^s \theta_i sign(2u-a_i-a_{i+1}) sign(2v-a_i-a_{i+1})\delta_i. \label{ideG}
\end{eqnarray}
Formula \eqref{ideG} defines the density a copula for any set $\{a_i, i=1,\cdots, s: 0\le a_1<a_2<\cdots<a_s<a_{s+1}= 1\}$. 
Via simple computations, we show that functions \[ f_i(x)=\frac{sign(2x-a_i-a_{i+1})}{\sqrt{a_{i+1}-a_i}} \] satisfy the conditions of Theorem \ref{T1} . Note that $\theta_i (a_{i+1}-a_i) =\lambda_i$ implies $\lvert \lambda_i \lvert \le 1/(a_{i+1}-a_i)$. Moreover, this copula coincides with copula \eqref{cops2} when $s=1$ and $a_1=0$. This copula departs from independence by adding local perturbations to rectangles $[a_i, a_{i+1}]\times[a_i, a_{i+1}]$ of size $\pm\theta$.

An interesting case is when all functions have the same maximum, same minimum and the absolute values of these numbers are same. The range of the $\lambda_k$ is symmetric. This is the case when $\varphi_k$ is $a\sin w_kx$, or $a\cos w_kx$. We consider these cases in the subsections below.

\subsubsection{Copulas based trigonometric on functions}
Here, we look into copula families based on trigonometric bases of $L^2(0,1)$. These copulas are different from those introduced in Chesneau (2021)\cite{Chesneau}, because they don't involve sums of variables. Moreover, the mixing structure of copulas of Chesneau (2021)\cite{Chesneau} is not as evident as that of the copula families that we construct here. It is also worth mentioning that the copula families that are introduced here can be finite or infinite sums and that their fold products remain in the considered classes. When modeling with these copulas, the existence of waves in the scatter plot can indicate the number of functions to consider in the sum. An observation on $\rho_S(C)$ and $\tau(C)$ shows that some terms can be added to increase linear correlation without affecting association.
\begin{itemize}
\item \textbf{The sine-cosine copulas as perturbations of $\Pi(u,v)$}
\end{itemize}
We call sine-cosine copulas perturbations of the independence copula that consider the trigonometric orthonormal basis of $L^2(0,1)$ that consists of functions: $\{1, \sqrt{2}\cos 2\pi kx, \sqrt{2}\sin 2\pi kx, k\in\mathbb{N^*}\}$. Theorem \ref{T1} guarantees that 
\begin{eqnarray} \label{copt}
C(u,v)=uv+\frac{1}{2\pi^2}\sum_{k=1}^\infty \frac{\lambda_k}{k^2} \sin 2\pi ku \times \sin 2\pi kv+\nonumber\\ 
+\frac{1}{2\pi^2}\sum_{k=1}^\infty \frac{\mu_k}{k^2}(1-\cos 2\pi ku)(1-\cos 2\pi kv)
\end{eqnarray}
is a copula, when $\displaystyle \sum_{k=1}^{\infty} (\lvert{}\lambda_k\lvert{}+\lvert{}\mu_k\lvert{})\le 1/2$. For this example, condition \eqref{lambda} holds automatically. The series \eqref{copt} also converges uniformly and absolutely on $[0,1]^2$ thanks to Weierstrass' $M$-test and the fact that the functions in the series are uniformly bounded on $[0,1]^2$.

\begin{figure}[h!]
\centering 
\includegraphics[scale=0.4]{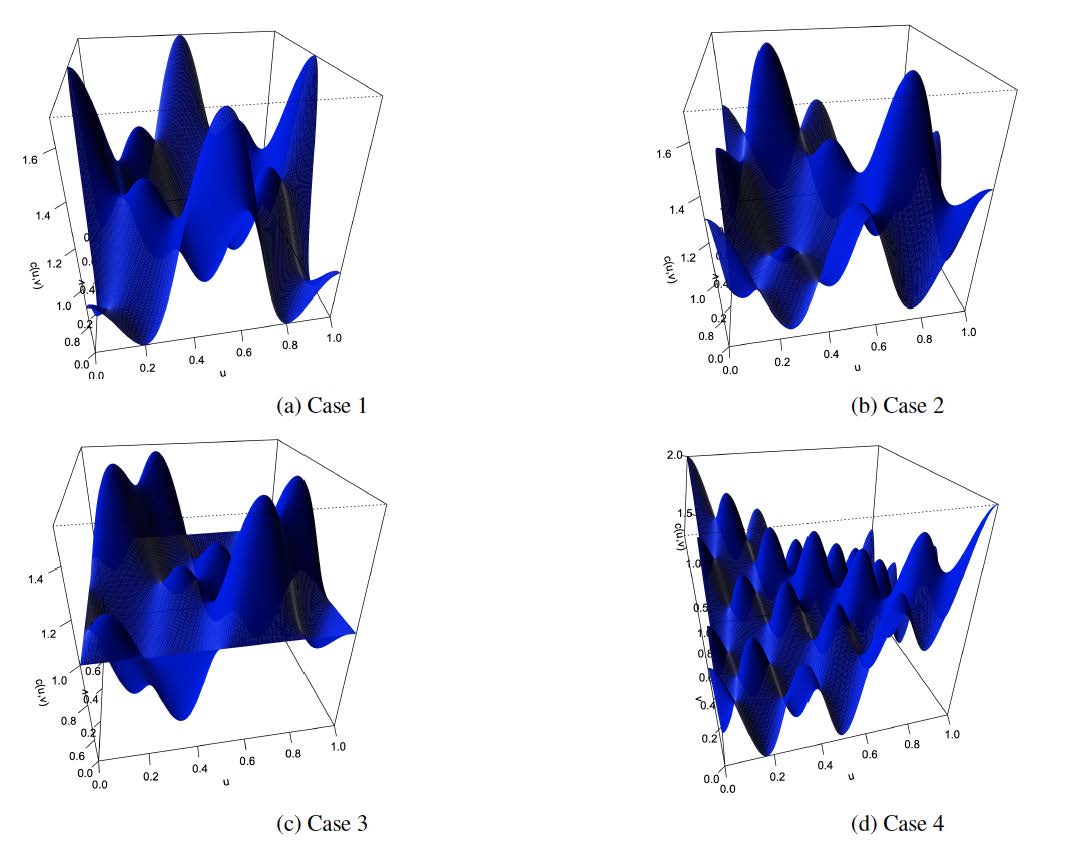}
\caption{Examples of sine-cosine copula densities}\label{sc}
\end{figure}

\noindent Figure \ref{sc}(a) represents the copula density \[ c(u,v)=1+.3\cos(4\pi u)\cos(4\pi v)-.4\cos(2\pi u)\cos(2\pi v), \] Figure \ref{sc}(b) is the graph of \[ c(u,v)=1+.3\cos(4\pi u)\cos(4\pi v)-.4\sin(2\pi u)\sin(2\pi v),\] Figure \ref{sc}(c) represents \[ c(u,v)=1+.3\sin(4\pi u)\sin(4\pi v)-.4\sin(2\pi u)\sin(2\pi v) \]
and Figure \ref{sc}(d) represents the copula density \[ c(u,v)=1+.2\cos(4\pi u)\cos(4\pi v)-.4\cos(2\pi u)\cos(2\pi v)\]\[+.4\cos6\pi u)\cos(6\pi v). \]

 The Spearman's $\rho$ and Kendall's $\tau$ for this copula family have been computed and are given below.
\begin{lemma} \label{sincosl} For any bivariate population with copula $C(u,v)$ defined by \eqref{copt}, the Spearman's $\rho_S(C)$ and Kendall's $\tau(C)$ are 
\begin{eqnarray} \nonumber
\rho_S(C)=\frac{6}{\pi^2}\sum_{k=1}^{\infty}\frac{\mu_k}{k^2}, \quad -\frac{3}{\pi^2}\le \rho_S(C)\leq \frac{3}{\pi^2}\\ 
\tau(C)=\frac{1}{\pi^2}\sum_{k=1}^{\infty}\frac{4\mu_k+2\lambda_k\mu_k}{k^2}, \quad -\frac{2}{\pi^2}\le \tau(C)\leq \frac{2}{\pi^2}.\nonumber
\end{eqnarray}
\end{lemma}

We can see on Figure \ref{sc} that the range of the copula decreases as we increase the number of terms in the sum. For sine-cosine copulas, Lemma \ref{sincosl} shows that terms related to cosine in the density of the copula don't affect $\rho_S(C)$, while terms with sine linearly affect $\rho_S(C)$. 

\begin{remark}
Any copula from this family with $\mu_k=0$ for all integer $k\ge 1$ satisfies $\rho_S(C)=\tau(C)=0$. This is a subclass of copulas of the form
\begin{equation}\label{subC}
C(u,v)=uv+\frac{1}{2\pi^2}\sum_{k=1}^\infty \frac{\lambda_k}{k^2} \sin 2\pi ku \times \sin 2\pi kv.
\end{equation}
Therefore, the Spearman's $\rho$ and Kendall's $\tau$ are not good parameters for tests of independence  because they don't uniquely identify a member of this family.
\end{remark}

\begin{itemize}
\item \textbf{The sine copulas as perturbations of $\Pi(u,v)$}
\end{itemize}

Consider the basis of $L^2(0,1)$ given by $\{1, \sqrt{2}\cos k\pi x, k>0,k\in\mathbb{N} \}$. By Theorem \ref{T1}, for any sequence $\{\lambda_k,k\in\mathbb{N}\}$ satisfying \eqref{cosGcond}, the function \eqref{cosG} is a copula.
\begin{equation}
C(u,v)=uv+\frac{2}{\pi^2}\sum_{k=1}^{\infty} \frac{\lambda_k}{k^2} \sin k\pi u \sin k\pi v \quad \text{with} \label{cosG}
\end{equation} 
\begin{equation}
\sum_{k=1}^\infty \lvert{}\lambda_k\lvert{}\leq 1/2. \label{cosGcond}
\end{equation}
\begin{figure}[h!]
\centering 
\includegraphics[scale=.35]{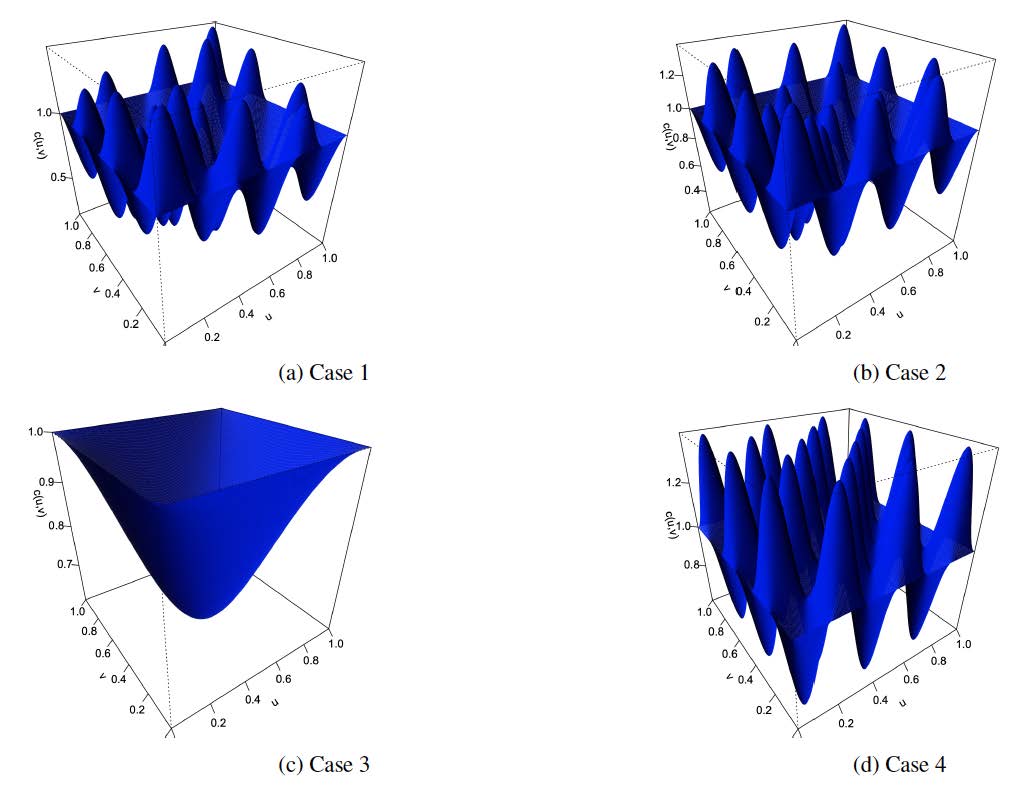}
\caption{Examples of sine copula densities}\label{sc1} 
\end{figure}

\begin{lemma}\label{lcop}
For copulas defined by formula \eqref{cosG}, we have
\begin{eqnarray}
\tau(C)=\frac{64}{\pi^4}\sum_{k=0}^{\infty}\frac{\lambda_{2k+1}}{(2k+1)^4}+\frac{16}{\pi^4}\sum_{k=1}^\infty\sum_{j=k+1}^\infty \frac{2\lambda_k\lambda_{j}}{(k^2-j^2)^2}(1-cos\pi(j+k)\pi)^2,\nonumber \end{eqnarray}
\begin{eqnarray}
\rho_S(C)=\frac{96}{\pi^4}\sum_{k=0}^{\infty}\frac{\lambda_{2k+1}}{(2k+1)^4}, \quad -\frac{48}{\pi^4}\leq \rho_S(C)\le \frac{48}{\pi^4}, \quad -\frac{32}{\pi^4}\leq \tau(C)\le \frac{32}{\pi^4}.\nonumber\end{eqnarray}
\end{lemma}
This family of copulas seems close to the sine-cosine family, but the formula of $\rho_S(C)$ or $\tau(C)$ in Lemma \ref{lcop} shows that they are not comparable. We see that for these copulas $\rho_S(C)=0$ when all odd coefficients are equal to zero. Formula \eqref{subC} defines a subclass of copulas from \eqref{cosG}, for which all odd coefficients are equal to zero.  Any finite sum of terms from \eqref{cosG} that includes the first term $uv$ and with coefficients that satisfy condition \eqref{cosGcond} is a copula. 
An example is 
\begin{equation}
C(u,v)= uv+\frac{2\lambda}{\pi^2}\sin\pi u \sin \pi v, \quad \lvert{}\lambda\lvert{}\leq 1/2. \label{cos}
\end{equation}

\begin{remark}
  The maximal coefficient of correlation of copula-based Markov chains generated by copulas \eqref{cos} is $\lvert{}\lambda\lvert{}$. We find $\rho_S(C)=\frac{96}{\pi^4}\lambda$ and $\tau(C)=\frac{64}{\pi^4}\lambda$.  So, $\lvert{}\tau(C)\lvert{}\le \lvert{}\rho_S(C)\lvert{}\le \lvert{}\lambda\lvert{}$. 
\end{remark}

\subsubsection{The extended Farlie - Gumbel - Morgenstern copulas}
In this subsection, we construct an example  of copula family based on formula \eqref{cop} that extends the Farlie-Gumbel-Morgenstern copulas. The FGM copula is a pertubation of the independence copula (see Longla et al (2022c) ). The form of the FGM copula allows an extension based on
 shifted Legendre polynomials, orthonormal basis of $L^2(0,1)$. It is important to mention that formula \eqref{cop} allows various families of extensions of the FGM copula that depend on the used orthonormal systems of functions. Legendre polynomials are defined by Rodrigues' formula 
\begin{equation}
P_k(x)=\frac{1}{2^k k!}\frac{d^k}{dx^k}(x^2-1)^k.
\end{equation}
This formula was discovered by Rodrigues in (1816) \cite{Ro} and leads to
\begin{equation}
P_k(x)=\frac{1}{2^k}\sum_{i=1}^{ [\frac{k}{2}] } (-1)^i {k \choose i}{2k-2i \choose k}x^{k-2i}.
\end{equation}
Legendre polynomials form a basis of $L^2(-1,1)$. When they are shifted to $[0,1]$, a renormalization is used to obtain shifted orthonormal polynomials. They appear in many fields of mathematics as solutions to a special differential equation and are used for approximation of functions. For these polynomials, $\int_{-1}^1P^2_k(x)dx=\frac{2}{2k+1}$. The transformation $y=2x-1$ gives to the orthonormal basis of shifted Legendre polynomials  as
\begin{equation}
\varphi_k(x)=\sqrt{2k+1}P_k(2x-1). \label{legp}
\end{equation}

For $k>0,$ $\varphi_k(1)=\sqrt{2k+1}=\max\varphi_k$, for odd $k$, $\min \varphi_k=-\sqrt{2k+1}=\varphi_k(0)$ and for even $k,$ $0>\min\varphi_k>-\sqrt{2k+1}$. This minimum is achieved at least at one point. Therefore, when \eqref{cop} contains $\varphi_k(x)$ alone, the range of $\lambda_k$ is \[ \frac{-1}{2k+1}\le \lambda_k\leq \min(1, \frac{1}{(2k+1)\lvert{}\min P_k(x)\lvert{}}). \] 

Based on the special structure of these functions, using the following equality known as Bonnet's recursion formula \[ (2n+1)P_n(x)=\frac{d}{dx}(P_{n+1}(x)-P_{n-1}(x)),\]
we establish the equality \[\displaystyle (2n+1)\int_{-1}^x P_n(t)dt= P_{n+1}(x)-P_{n-1}(x).\] This formula uses the fact that $P_n(-1)=(-1)^n$. Substituting $u=\frac{1}{2}(x+1)$ and $t=2s-1, 2ds =dt$, we obtain
\[ \displaystyle (2n+1)\int_{0}^{u} P_n(2s-1)2ds= P_{n+1}(2u-1)-P_{n-1}(2u-1). \]
Thus, using the definition of shifted Legendre polynomials $\varphi_k(x)$, we have
\begin{eqnarray}
2\sqrt{2n+1}\int_{0}^{u} \varphi_n(s)ds= \frac{1}{\sqrt{2n+3}}\varphi_{n+1}(u)-\frac{1}{\sqrt{2n-1}}\varphi_{n-1}(u). \label{lege}
\end{eqnarray}
For more on Legendre polynomials, see Belousov (1962) \cite{Belousov}, Sz\"ego (1975) \cite{Szego} and the references therein. Formula \eqref{lege} is used to establish the following. 
\begin{lemma} \label{LFGM}
For any copula of the extended FGM copula family, 
\[ \rho_S(C)= \lambda_1 ,\quad \text{and} \quad \tau(C)=\frac{2\lambda_1}{3}-\sum_{k=1}^\infty \frac{2\lambda_k\lambda_{k+1}}{(2k+1)\sqrt{2k+3}\sqrt{2k-1}}. \]
\end{lemma}
Simple computations show that $-1/3\le \rho_S(C)\le 1/3$ and $\tau(C)$ can be increased or decreased by adding more non-zero terms to the series. In general, copulas with $\varphi_k(x)$ given in formula \eqref{legp} satisfy condition \eqref{cond2}. We call this set of copulas extended Legendre-Farlie-Gumbel-Morgenstern family of copulas. Some copulas of this form are as follows. For $\lvert{}\lambda\lvert{}\le 1/3, $ $3\lvert{}\lambda_1\lvert{}+5\lvert{}\lambda_2\lvert{}\le 1$,
\begin{eqnarray} 
C_{\lambda}(u,v)=uv+3\lambda (u-u^2)(v-v^2), \label{FGM1}\\ 
C_{\lambda_1, \lambda_2}(u,v)=uv+3\lambda_1 (u-u^2)(v-v^2)+\nonumber \\ +5\lambda_2(2u^3-3u^2+u)(2v^3-3v^2+v). \label{legex}
\end{eqnarray}

The first example above is the FGM copula. For more on this copula family, see Farlie (1960) \cite{FARLIE}, Gumbel (1960) \cite{GUMBEL}, Jonhson and Kotz (1975) \cite{Johnson} or Morgenstern (1956) \cite{Morgenstern}. In general this copula has parameter $\theta$ in place of $3\lambda_1$. The second example is one of the possible extensions of this copula family. For this extension, we have 
\begin{equation}\label{dd}
\rho_S(C_{\lambda_1,\lambda_2})=\lambda_1, \quad \tau(C_{\lambda_1,\lambda_2})=\frac{2}{3}\lambda_1(1-\frac{\lambda_2}{\sqrt{5}}).
\end{equation}

Any subclass of the copula family defined with Legendre polynomials and containing the copula given by formula \eqref{FGM1} is an extension of the FGM family. Depending on what the investigator is looking for, it might be better to set some coefficients equal to zero, deal with the signs of the coefficients to increase or reduce the strength of dependence. This is justified by formula \eqref{dd}.  The FGM copula family has parameter $\theta=3\lambda$, with $\lvert \theta \lvert \le 1$; so $\lvert \lambda \lvert \le 1/3$. No copula from our extension can have $\rho_S(C)$ out of this range because increasing the number of non-zero parameters reduces the range of $\lambda_1$. Example \eqref{legex} has $-1/5\le \lambda_2\leq 1/5$ when $\lambda_1=0$, but can be extended to $-1/5\le \lambda_2\leq 2/5$. 

\begin{figure}[ht!]
\centering
\includegraphics[scale=.32]{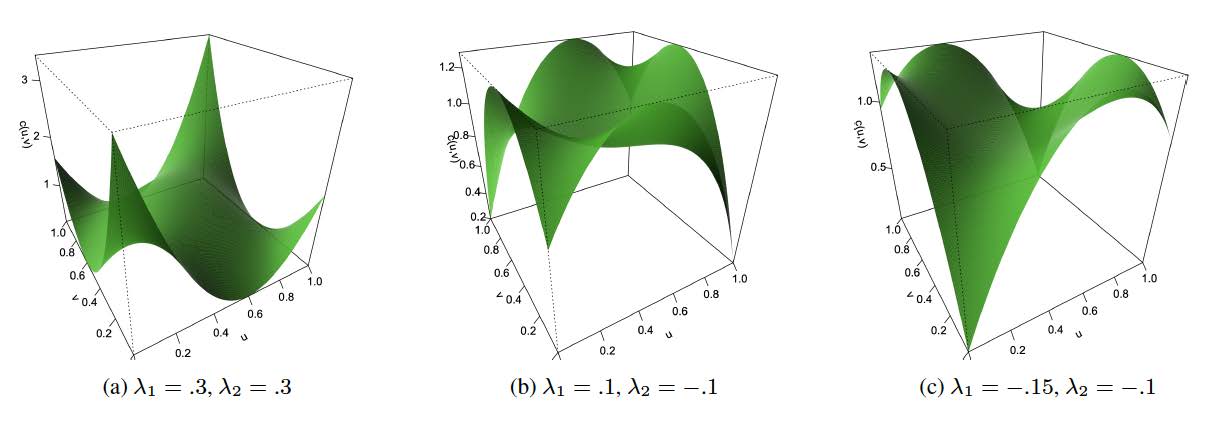}
\caption{Densities of Legendre-Farlie-Gumbel-Morgenstern copulas}
\end{figure}
For any copula from this family, it is worth formulating the following.
\begin{remark} Lemma \ref{LFGM} shows that any copulas from this family with the same $\lambda_1$ have the same $\rho_S(C)$, while it is not the case for $\tau(C)$. Other formulations of these copulas  use the basis $\varphi_{k}(x)=\sqrt{2j_k+1}P_{j_k}(2x-1)$, where $j_k$ is a sequence of positive integers.
\end{remark}

\subsubsection{Counter example - the sine basis}

This case is to show that not all orthonormal bases of $L^2(0,1)$ can be used to construct copulas with square integrable densities. We take the case of functions constructed using the sine basis. They are not copulas under any assumptions. In fact, if
\begin{equation}
C(u,v)=\frac{2}{\pi^2}\sum_{k=1}^{\infty} \frac{\lambda_k}{k^2} (1-\cos k\pi u)(1- \cos k\pi v), \quad \text{then}
\end{equation} 
\begin{equation}
c(u,v)=2\sum_{k=1}^{\infty} \lambda_k \sin k\pi u \sin k\pi v.
\end{equation}
If this function is a copula, then using $C(u,1)=u$ and taking the derivative of the series, we obtain
\begin{equation}
1=2\sum_{k=1}^{\infty} \frac{\lambda_k}{k\pi} \sin k\pi u (1-\cos k\pi ).
\end{equation}
Multiplying by $\sin k\pi u$ and integrating leads to $\frac{1-\cos k\pi}{k\pi} =\frac{\lambda_k}{\pi k} (1-\cos k\pi).$ Therefore, $\lambda_{2k+1}=1$. As said earlier, the function $\varphi(x)=1$ in this case is embedded in a subspace of dimension $\infty$ of eigenfunctions associated to the eigenvalue $1$. This is in contradiction with condition \eqref{lambda}. So, no such decomposition represents a square integrable copula density. 

\section{Copula-based Markov chains and $\psi$-mixing}
We study here the long run behavior of several copula families constructed with various sets of bases of $L^2(0,1)$. We provide the study of the mixing structure of copula-based Markov chains generated by our constructed copulas. In some cases, we provide a modification of the copulas to present them as convex combinations of some members of the constructed families. This transformation provides a link between this work and results of Longla et al. (2022) \cite{L0} and Longla (2015) \cite{Lo2015}.
Mixing coefficients are used in probability theory to establish central limit theorems for sums of dependent random variables (see Peligrad and Utev (1997)\cite{PelUt}), which are used to construct confidence intervals for means of functions of Markov chains. Copulas characterize the dependence structure of Markov chains or multivariate random variables in general via Sklar's theorem (See Sklar (1959) \cite{SK}) when marginal distributions are continuous. Therefore, it is important to look into some conditions on copulas that would guarantee a given dependence structure (see Beare (2010) \cite{Beare}, Longla (2015)\cite{Lo2015} among others). Dependence structures considered in the literature include $\alpha$-mixing, $\beta$-mixing, $\rho$-mixing, $\phi$-mixing, $\psi$-mixing, $\psi'$-mixing, $\psi^*$-mixing and others. We define only mixing coefficients of interest in this paper. For more on mixing coefficients, see Bradley (2007) \cite{BR}. 
For a stationary Markov chain $X_0, \cdots, X_n$, $\mathbb{A}=\sigma(X_0)$, mixing coefficients are defined as
\[
\phi_n=\sup_{P(A)P(B)>0, A,B\in\mathbb{A}} \lvert{}P^n(B\lvert{} A)-P(B)\lvert{}, \] \[\rho_n=\sup_{f}\{\mathbb{E}(f(X_0)f(X_{n})): \mathbb{E}f(X_0)=0, \mathbb{E}f^2(X_0)=1\},\]
\[
\psi_n=\sup_{P(A)P(B)>0, A,B\in\mathbb{A}}\frac{\lvert{}P^n(A\cap B)-P(A)P(B)\lvert{}}{P(A)P(B)},
\]
\[
\psi'_n=\inf_{P(A)P(B)>0, A,B\in\mathbb{A}}\frac{P^n(A\cap B)}{P(A)P(B)} \quad \text{and}\] \[ \psi_n^{*}=\sup_{P(A)P(B)>0, A,B\in\mathbb{A}}\frac{P^n(A\cap B)}{P(A)P(B)},
\] 
where $P^n(B\lvert{} A)$ stands for $P(X_n\in B\lvert X_0\in A)$, and $P^n(A\cap B)$ stands for $P(X_0\in A ~ and ~ X_n\in B)$.
 Using these coefficients, Bradley (2005) \cite{Bradley} in his survey presented the following result.
Suppose $X_0, \cdots, X_n$ is a strictly stationary sequence of random variables. If there exists $n$ such that $\psi'_n>0$ and $\psi_n^{*}<\infty$, then $\psi_n\to 0$ (the sequence is $\psi$-mixing). Based on this statement, results of Longla (2015) \cite{Lo2015} and Longla et al. (2022c) \cite{LoMous}, we can conclude the following.
\begin{theorem} \label{theoM}
An absolutely continuous copula C(u,v) generates $\psi$-mixing stationary Markov chains when there exist integers $s_1$ and $s_2$ such that some versions of the densities of $C^{s_1}(u,v)$ and $C^{s_2}(u,v)$ satisfy one of the following conditions.
\begin{enumerate}
\item There exists a constant real number $M$ such that $c^{s_1}(u,v)>0$ on a set of Lebesgue measure $1$ and $c^{s_2}(u,v)<M $ for all $(u,v)\in [0,1]^2$.
\item For all $(u,v)\in [0,1]^2$, $c^{s_2}(u,v)<M\le 2 $.
\end{enumerate}
\end{theorem}

It is important in Theorem \ref{theoM} that the version of the density is bounded above on $[0,1]$. Failure of this condition on a set of Lebesgue measure $0$ can imply that there is no $\psi$-mixing. 
\begin{example}\label{rem1}
The copula $C(u,v)=aM(u,v)+(1-a)W(u,v),$ with $0\le a\le 1,$ from the Mardia family, when used with continuous marginals, doesn't generate $\rho$-mixing Markov chains as shown in Longla (2014) \cite{L2}. These copulas are not absolutely continuous. Another copula that fails these conditions is \eqref{cops2} with $\alpha=1=\lambda$.
\end{example}
This example is a convex combination of $W(u,v)=\max(u+v-1,0)$ and $M(u,v)=\min(u,v)$. 
A convex combination of copulas $C_1(u,v)$, $\cdots,$ $C_k(u,v)$ is the copula
$
C(u,v)= a_1 C_1(u,v)+ \cdots a_k C_k(u,v), 
$
with $0< a_i,$ for all $i$ and $a_1+\cdots+a_k=1.$
Using Theorem \ref{theoM}, we get.
\begin{corollary}
Let $C_1(u,v), \cdots, C_k(u,v)$ be copulas such that the density of $C_{s_1}* C_{s_1+1}*\cdots*C_{s_2}(u,v)$ is bounded away from zero on a set of Lebesgue measure $1$ for some positive integers $s_1, s_2$. Assume that a version of the density of $(a_1C_{1}+\cdots+a_kC_{k})^{s_3}$ is bounded on $[0,1]^2$ for some positive integer $s_3$ or a version of the density of $(\sum a_iC_i)^{s_4}(u,v)$ is strictly less than $2$ for some $s_4\in \mathbb{N}$. Then, any convex combination of copulas $C_1(u,v), \cdots, C_k(u,v)$ generates exponential $\psi$-mixing stationary Markov chains with continuous marginals.
\end{corollary}

\subsection{Applications to our new copula families}
In this subsection, we apply the result on mixing to Markov chains generated by our newly created copula families. 
For any copula of the form \eqref{cops2}, the maximal coefficient of correlation is $\lvert \lambda\lvert$. Therefore, it generates $\rho$-mixing Markov chains when $\lvert{}\lambda\lvert{}<1$. Moreover,
\begin{enumerate}
\item For $\alpha\ne 1$, $C_{\alpha,\lambda}(u,v)$ generates $\psi$-mixing when $-\min(\frac{1}{\alpha}, \alpha) \le \lambda<1$;
\item $C_{1,\lambda}(u,v)$ generates $\psi$-mixing when $-1 < \lambda<1$;
\item There is no $\psi'$-mixing, $\psi^*$ or $\psi$-mixing for other values of $(\alpha,\lambda)$.
\end{enumerate}

Note that when $\alpha=1$, the range of $\lambda$ includes $\pm1$, for which there is no $\rho$-mixing, $\psi^*$-mixing or $\psi'$-mixing. This case is equivalent to

\begin{equation}
c_3(u,v)=1+\lambda \cdot sign(2u-1)\cdot sign(2v-1), \quad\quad \lvert{}\lambda\lvert{}\le 1. \label{ide}
\end{equation}
When $\lvert{}\lambda\lvert{}< 1$, the density \eqref{ide} splits $[0,1]^2$ into 4 subsets. It is constant on each of the subsets, and the union of the 4 subsets is of full Lebesgue measure. Clearly, for $\lvert{}\lambda\lvert{}=1$, this density is not bounded away from zero. 

For a stationary copula-based Markov chain $U_0, \cdots, U_n$ based on \eqref{cop} with uniform marginals, the copula of $(U_0,U_n)$ is $C^n(u,v)$ (see Longla et al. (2022)\cite{L0}). 
\begin{theorem}\label{mixcop}
For any copula defined by \eqref{cop} and satisfying condition \eqref{cond},
\begin{equation}
C^{n}(u,v)=uv+\sum_{k=1}^{\infty}\lambda_k^n \int_{0}^u \varphi(x)dx\int_0^v\varphi(y)dy. \label{joint}
\end{equation}
Moreover, $\rho_n=\sup_{k}\lvert{}\lambda_k\lvert{}^n$ and Markov chains generated with uniform marginals are both geometrically ergodic and exponential $\rho$-mixing if and only if $\sup_k \lvert{}\lambda_k\lvert{}<1$. If the sum has a finite number of non-zero terms with $\lvert \lambda_k\lvert<1$ for all $k$, then the Markov chains generated with uniform marginals are $\psi$-mixing. 
\end{theorem}
\begin{corollary}
For any set $\{a_i, i=1,\cdots, s: 0\le a_1<a_2<\cdots<a_s<a_{s+1}= 1\}$, and $\lvert{}\theta_i\lvert{}\leq 1$, formula  \eqref{ideG} defines the density of a copula that generates $\psi-mixing$ stationary Markov chains with continuous marginal distributions.
\end{corollary}

\subsection{Mixing properties of the extended FGM family}
We investigate the long run behavior of Markov chains generated by extended FGM copulas. Define $\alpha_k=1$ if $k$ is odd and $\alpha_k=\min P_k(x)$, if $k$ is even.

\begin{theorem}\label{decompPsi}
For large enough values of $n$, a version of the copula density $c^n(u,v)$ is bounded above by $M<2$ when $\varphi_k(x)$ is defined using Legendre polynomials with $\sum{\lambda_k^2}<\infty$, $\sup_k\lvert{}\lambda_k\lvert{}<1$ and $\sum\lvert{}\lambda_k\lvert{}(2k+1)\alpha_k\le 1$. Therefore the copula $C(u,v)$ generates $\psi^*$-mixing stationary Markov chains (which is equivalent to $\psi$-mixing for stationary Markov chains) for all values of its parameters. 
\end{theorem}

\begin{remark}
It follows from the proof of Theorem \ref{decompPsi}, that $c^{n}(x,y)\to 1$ as $n\to\infty$, when the conditions of Theorem \ref{decompPsi} are satisfied.
For any system of eigenfunctions $\varphi_k(x)$ and for all values of the parameters described in Theorem \ref{decompPsi}, Theorem \ref{mixcop} holds because $\psi$-mixing implies $\rho$-mixing. Theorem \ref{mixcop} holds when extended FGM copulas are considered. A result similar to Theorem \ref{decompPsi} is valid for our trigonometric extensions \eqref{copt} and \eqref{cosG} in the form: for $\sum(\lvert{}\lambda_k\lvert{}+\lvert{}\mu_k\lvert{})\le 1/2,$ generated Markov chains are $\psi$-mixing.
\end{remark}
Based on the analysis above, we can derive several facts about mixing properties of Markov chains generated by absolutely continuous copulas with square integrable densities. In general, if we assume that a Markov chain is generated by a copula \eqref{cop} satisfying conditions \eqref{lambda} and \eqref{cond2} with an absolutely continuous marginal distribution, then the following holds. 
\begin{theorem} \label{TT} Under the assumptions of Theorem \eqref{mixcop}, 
\begin{enumerate}
\item If the sequence $\varphi_k(x)$ is uniformly bounded, and $\sup_k\lvert{}\lambda_k\lvert{}<1$ for all $k$, then the Markov chain is $\psi$-mixing.
\item If $\sup_k\lvert{}\lambda_k\lvert{}<1$ for all $k$, then the Markov chain is $\rho$-mixing.
\item If Condition \eqref{cond2} allows strict inequality, then the Markov chain is $\psi'$-mixing.
\end{enumerate}
\end{theorem}
Each of the examples of Figure \eqref{sc} is bounded away from zero on $[0,1]^2$ and is strictly less than $2$. Therefore, by Theorem \eqref{TT}, these examples generate $\psi$-mixing Markov chains. This means that the sine copulas and sine-cosine copulas generate $\psi$-mixing Markov chains. This is justified by the fact that each of the basis functions is bounded by $\sqrt{2}$ and condition \eqref{cond2} implies $\lvert \lambda_k\lvert\le 1/2$.
It is important to note that the uniform bound on the sequence $\varphi_k(x)$ is crucial in statement 1 of Theorem \ref{TT}. 

\begin{figure}[h!]
\centering
\includegraphics[scale=.32]{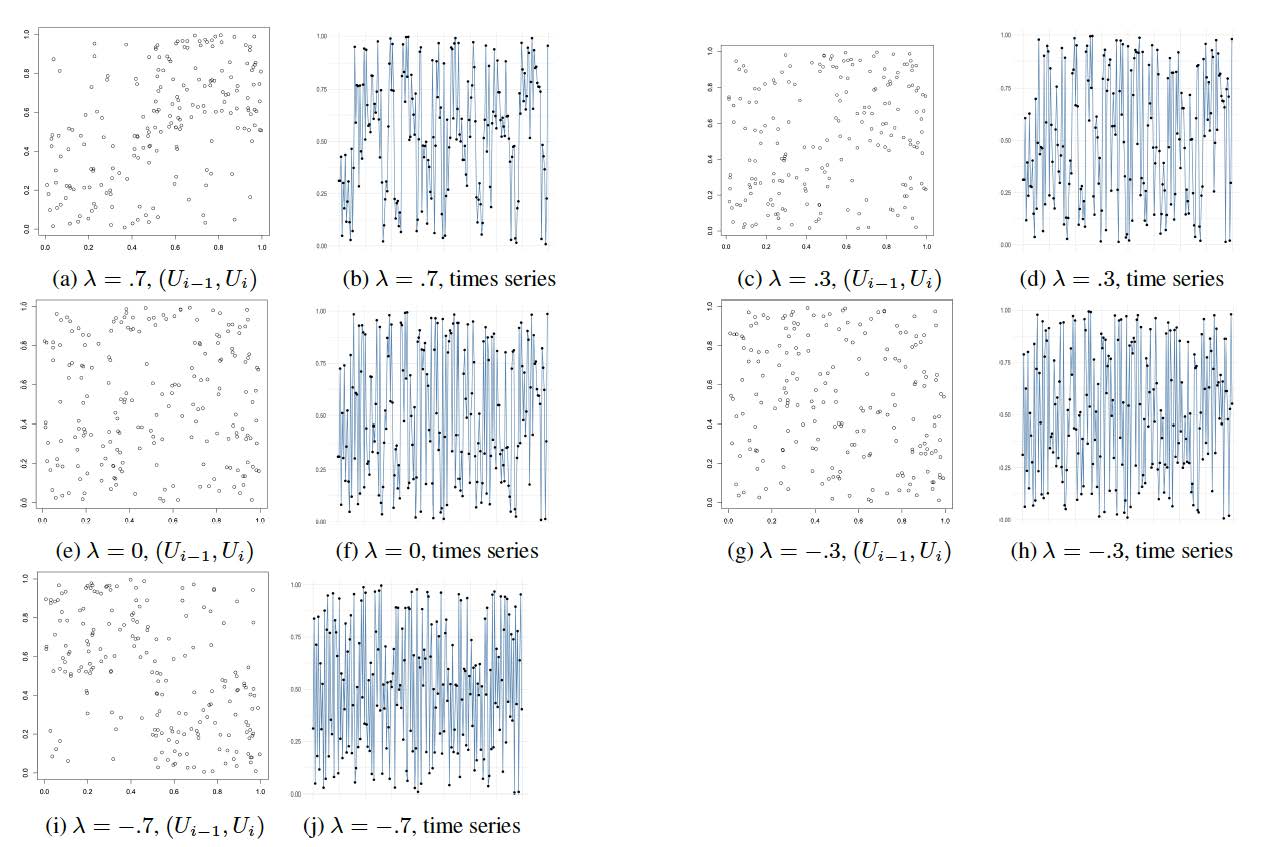}
\caption{Density of copula $C_{1,\lambda}$}
\end{figure}
\subsection{Answer to an open question on mixing }
Longla (2015) was solely devoted to mixing properties of copula-based Markov chains. Results were provided for several mixing coefficients, in one direction. There was a proof that mixing for Markov chains generated by copulas implies mixing for Markov chains generated by their convex combinations. In this subsection, we show that mixing for Markov chains generated by convex combinations doesn't require mixing for Markov chains generated by any of the copulas of the convex combination. To do this, we take extreme cases of formula \eqref{ide}, defined by $\lambda=\pm1$. None of them exhibits mixing in the sense of the coefficients defined in this paper, but they have square integrable densities. 
\begin{remark} We provide here a connection between various results of Longla et al. (2022) \cite{L0}, Longla (2015) \cite{Lo2015} and Longla et al (2022c) \cite{LoMous} on mixing.
\begin{enumerate}
\item The definition of $C_3(u,v)$ through its density leads to 
\[ C_3(u,v)=\frac{1+\lambda}{2}C_1(u,v)+\frac{1-\lambda}{2}C_2(u,v), \]
where $C_1(u,v)$ and $C_2(u,v)$ are the copulas of this family defined by $\lambda=\pm1$. Moreover, $C_3(u,v)$ generates $\psi$-mixing Markov chains for $\lvert{}\lambda\lvert{}<1$. 
\item None of the two copulas $C_1(u,v)$ and $C_2(u,v)$ has a density bounded away from zero on a set of Lebesgue measure 1. Moreover, $C_1(u,v)$ and $C_2(u,v)$ don't generate $\psi'$-mixing Markov chains with continuous marginal distributions. This is because they don't generate $\rho$-mixing Markov chains with continuous marginals. Their convex combinations with $\lvert{}\lambda\lvert{}< 1$ generate $\psi'$-mixing (which implies that they generate $\rho$-mixing).
\item None of the copulas $C_1(u,v)$ and $C_2(u,v)$ generates $\psi$-mixing, because they would otherwise generate $\rho$-mixing; but their convex combinations with $\lvert{}\lambda\lvert{}< 1$ generate $\psi$-mixing Markov chains.
\item None of the copulas $C_1(u,v)$ and $C_2(u,v)$ generates $\phi$-mixing, because they would otherwise generate $\rho$-mixing; but their convex combinations with $\lvert{}\lambda\lvert{}<1$ generate $\phi$-mixing Markov chains.
\item None of the copulas $C_1(u,v)$ and $C_2(u,v)$ generates $\beta$-mixing, because they would otherwise generate $\rho$-mixing; but their convex combinations with $\lvert{}\lambda\lvert{}< 1$ generate $\beta$-mixing Markov chains. This is because they are symmetric and generate $\rho$-mixing (see Longla and Peligrad (2012)) \cite{LoMag}.
\end{enumerate}
\end{remark}
It is crucial in the above results that the copulas are absolutely continuous. To see this, take the copula
\[C_{\alpha,\theta}(u,v)=\alpha(uv+\theta (u-u^2)(v-v^2))+(1-\alpha) \min(u,v)\] \[ =\alpha uv + \alpha \theta (u-u^2)(v-v^2)+(1-\alpha)\min(u,v).\]
$C_{\alpha,\theta}(u,v)$ is a convex combination of copulas containing $\Pi(u,v)$. Therefore, by Longla (2015) \cite{Lo2015}, it generates $\psi{'}$-mixing stationary Markov chains. $C_{\alpha,\theta}(u,v)$ is a convex combination that contains $M(u,v)$. Therefore, by Longla et al. (2022) \cite{L0}, it doesn't generate $\psi$-mixing stationary Markov chains with continuous marginals. Its density $c_{\alpha,\theta}(u,v)=\alpha(1+\theta(1-2u)(1-2v))<2$ for all $(u,v)\in [0,1]^2$ et all $\alpha$ and $\theta$, such that $0<\alpha<1, 0<\lvert{}\theta\lvert{}<1$. Therefore, the conditions cannot be weakened. 

Formula \eqref{cops2} defines an absolutely continuous copula that generates $\psi$-mixing stationary Markov chains. Note that its density is equal to 0 on a set of non-zero Lebesgue measure for $\lambda=-\alpha>-1$ ($\lambda<1$)or $\lambda=-\frac{1}{\alpha}>-1$ ($\lambda>1$). This by itself answers an open question on $\psi$-mixing from Longla (2015)\cite{Lo2015}. Longla (2015) showed that when the density is bounded away from $0$ on a set of Lebesgue measure $1$, the copula generates $\psi'$-mixing. This example shows that the condition is very strong. We can obtain $\psi$-mixing even for copulas with density equal to zero on a set of strictly positive Lebesgue measure using copulas from these new families.  

\section{Central limit theorem and simulation study}
In this section we consider some functions of the Markov chain and derive central limit theorems for estimators of parameters of the copula and the population mean. We use an example of copula from the derived copula families, with special association properties. In fact, we select a copula that has $\rho_S(C)=0$, and a modification that also includes $\tau(C)=0$. For these examples, regular estimation procedures based on these measures of association fail. 
\subsection{Central limit theorem for parameter estimators}
We consider in this section central limit theorems for averages of functions of the Markov chain generated by the copula 
\begin{equation}
C(u,v)=uv+\mu_1 \Phi_1(u)\Phi_1(v)+\mu_2 \Phi_2(u)\Phi_2(v) \label{model},
\end{equation} 
with $\varphi_1(x)=\sqrt{2}\sin 2\pi x$, $\varphi_2(x)=\sqrt{2}\sin 4\pi x$ and $2\lvert \mu_1\lvert+2\lvert \mu_2\lvert \leq 1$. Let $f(x)$ be such that $\mu=\int_0^1f(x)dx$ and $\int_0^1f^2(x)dx<\infty$. Define $S_n(f)=\sum_{i=1}^n f(U_i)$, $$\hat{\mu}_k=\frac{1}{n-1}\sum_{i=1}^{n-1}\varphi_k(U_i)\varphi_k(U_{i+1}), \quad \text{for}\quad k=1,2.$$
\begin{theorem}Assume that $U_1, \cdots U_n$ is a realization of the stationary Markov chain generated by \eqref{model} and the uniform marginal distribution. The following holds.
\begin{enumerate}
\item \begin{equation}  
\sqrt{n-1}\large\left( {\hat{\mu}_1\choose \hat{\mu}_2}-{\mu_1\choose \mu_2}\right)\to N\left( {0 \choose 0} , \begin{pmatrix}1 & -\mu_1\mu_2 \\ -\mu_1\mu_2 & 1\end{pmatrix}\large\right).
\end{equation}
\item  The central limit theorem holds in the form 
\begin{equation}
\sqrt{n}(\frac{S_n(f)}{n}-\mu) \to N(0, \sigma^2_f), \quad \text{where},\quad \sigma^2_f=\sigma^2+2(\frac{\mu_1 A_1^2}{1-\mu_1}+\frac{\mu_2A_2^2}{1-\mu_2}),
\end{equation}
where $A_i^2=2(\int_0^1\sin 2i \pi u f(u)du)^2$.
\end{enumerate} \label{mut}
\end{theorem}
Based on properties of the copulas that we have constructed here, the variance $\sigma^2_f$ exists, is finite and is strictly positive.
To construct a confidence set or test a specific value of $(\mu_1,\mu_2)$, we can use the following consequence of Theorem \ref{mut}.

\begin{corollary}\label{colmu}
Assume that $U_1, \cdots U_n$ is a realization of the stationary Markov chain generated by \eqref{model} and the uniform marginal distribution. The following holds. 
\begin{equation}\nonumber
X_n=\frac{(n-1)}{1-{\mu}_1{\mu}_2}\large\left( (\hat{\mu}_1-\mu_1)^2+(\hat{\mu}_2-\mu_2)^2+2{\mu}_1{\mu}_2(\hat{\mu}_1-\mu_1)(\hat{\mu}_2-\mu_2)\large\right)
\end{equation}
has an approximately $\chi^2(2)$ distribution. For confidence interval purposes, replace $\mu_1\mu_2$ by $\hat{\mu}_1\hat{\mu}_2$.
\end{corollary}
Consider now $\lambda_k=0$ for all $k$ and $\mu_k=0$ for $k>2$. Moreover, request that $\tau(C)=0$ for the sine-cosine copula. This implies $\mu_2=-4 \mu_1$, reducing the copula to the one parameter family with $\lvert \mu_1\lvert\le 0.11$ 
\begin{equation}
C(u,v)=uv+\frac{\mu_1}{2\pi^2}[(1-\cos 2\pi u)(1-\cos 2\pi v)-(1-\cos 4\pi u)(1-\cos 4\pi v)]. \label{taurho}
\end{equation}
 \begin{remark} Consider a Markov chain generated by \eqref{taurho}. Theorem \ref{mut} implies the following.
\begin{enumerate}
\item For $f(x)=\mathbb{I}(x\le a)$, $a\in(0,1)$, we get a dependent Bernoulli sequence and a CLT for the sample relative frequency of success holds with \begin{equation}
 \sigma^2_f= a(1-a)+4\frac{\mu_1}{\pi^2}\large\left(\frac{\sin^4 \pi a}{1-\mu_1}-\frac{\sin^4 2\pi a}{1+4\mu_1}\large\right). \label{I}
\end{equation}
Note that for $\mu_1\ne 0$, the observations are dependent, but there exists a value of $a=\frac{1}{\pi}\arccos((\frac{1-4\mu_1}{4(1-\mu_1)})^{1/4})$ for which $\sigma^2_f=a(1-a)$. This means that the sample behaves as if it was a simple random sample. A researcher who starts with this simple random sample assumption would have wrong conclusions. 
\item\noindent For $f(x)=-\lambda\ln(1-x)$, we get a Markov chain with exponential marginal distribution of parameter $\lambda$. The estimator of $\lambda$ is $S_n(f)/n$ and the CLT holds with  $(\int_0^1\sin 2\pi t \ln (1-t)dt)^2=0.1505165$, $4(\int_0^1\sin 4\pi t \ln (1-t)dt)^2=0.245684$ (these values where obtained using $R$ to integrate). \begin{equation}
\sigma^2_f=\lambda^2+4\mu_1\lambda^2\large\left(\frac{0.1505165}{1-\mu_1}-\frac{0.245684}{1+4\mu_1}\large\right). \label{I2}
\end{equation} 
This example shows that we can still estimate $\lambda$ using the sample mean, but the variance of the estimator has to be treated more carefully. The asymptotic variance is not $\lambda^2/n$. It is multiplied by a coefficient that can be less than or grater than 1 depending on the choice of $\mu_1$.
\item For $f(x)=x$, the asymptotic variance is always larger than it would be, if the data was assumed from a simple random sample.   \begin{equation}
\sigma^2_f=\frac{1}{12}+\frac{\mu_1/\pi^2}{1-\mu_1}-\frac{\mu_1/\pi^2}{1+4\mu_1}=\frac{1}{12}+\frac{5\mu_1^2}{\pi^2(1-\mu_1)(1+4\mu_1)}.\label{vm}
\end{equation}
\item For $ f(x,y)=2w \sin(2\pi x)\sin(2\pi y)-\frac{1}{2}(1-w)\sin(4\pi x)\sin(4\pi y)$,  $w\in[0,1]$, we have the estimator 
$$\mu_w=\sum_{i=2}^{n}\large\left(2w \sin(2\pi U_{i-1})\sin(2\pi U_i)-\frac{1}{2}(1-w)\sin(4\pi U_{i-1})\sin(4\pi U_i)\large\right).$$
The central limit theorem for this estimator holds with variance less than that of $\hat{\mu}_1$,
\begin{equation}
\sigma^2_f= 1-2(1-4\mu_1^2)(w-w^2).
\end{equation} 
The smallest possible value of this variance is $\sigma^2_0=0.5+2\mu_1^2$. 
\end{enumerate}
\end{remark}
\subsection{Simulation study}
We generate a Markov chain ($U_1,\cdots U_{1000}$) using the uniform distribution on (0,1) and the copula \eqref{taurho} for $\mu_1=0.05$.
We use the following standard procedure to generate the Markov chain. 
Randomly generate an observation from $Uniform(0,1)$ as $U_1$. and for every $i>1$, we generate a new observation $W$ from $Uniform(0,1)$ and solve for $U_{i}$ the equation $C_{,1}(U_{i-1}, U_i)=W$. This step is implemented using the R function $uniroot$. After generating this sample, we performed a correlation test. This test showed that one could conclude there is no correlation between variables along the Markov chain. Moreover, the scatter plot would also convince further the investigator that this data seems to be a random sample. We use $\alpha=0.05$ for level of confidence of the intervals.
\begin{enumerate}
\item Case 1: We use $a=01., 0.2, 0.3, 0.4, 0.5, 0.6, 0.7, 0.8, 0.9$ to obtain a set of variables  $Yi=\mathbb{I}(U\le 0.i)$ for our sample. For this function of the observations we get a sample of values $(Y1\dots Y9)$. We now consider the problem of estimation of the parameters $a_i$ of the model using the obtained dependent Bernoulli sequence. We get $\hat{a}_i=\bar{Y}i$ and the variance is given by formula \eqref{I} for every $i$. A confidence interval is constructed. We run $N$ such studies, and get the proportions of intervals that contain the true value of the parameter. This is done under independence assumption and under our model. We obtain the following comparative results, that report the coverage probabilities under each of the models and for the given values of $100$ and $1000$ for $N$.

\begin{table}[h!]
\centering
\captionsetup{justification=centering,margin=0.5cm}
\begin{tabular}{|c||c|c|c|c |c|c|c|c|c| }
\hline
N/a &0.1 & 0.2 & 0.3 & 0.4 & 0.5 & 0.6 & 0.7 & 0.8 & 0.9\\ 
\hline\hline
100 & 95 &90 &94 &91 &91 &93 &96&97 &96\\ 
\hline 
100 & 91& 99& 95& 94& 94& 97 &98 &96 &96\\
\hline
1000 &95.7& 95.0& 94.8& 93.3& 92.9 & 92.7 & 94.5 &95.0 &95.5\\
\hline
\end{tabular}
\caption{Coverage probabilities under assumption of independence} 
\end{table}

\begin{table}[h!]
\centering
\begin{tabular}{|c||c|c|c|c |c|c|c|c|c| }
\hline
N/a &0.1 & 0.2 & 0.3 & 0.4 & 0.5 & 0.6 & 0.7 & 0.8 & 0.9\\ 
\hline\hline
100 & 92 &95 &96 &96 &96 &95 &94&94 &93\\ 
\hline 
100 & 96& 92& 95& 95& 93& 93& 91& 96& 96\\   
\hline 
1000 & 94.6& 95.0& 95.0& 96.4& 95.3& 93.4 &92.9 &94.3 &94.7\\
\hline
\end{tabular}
\caption{Coverage probabilities under our model C(u,v)}
\end{table}
 These tables show that the researcher would be off by a lot by assuming independence for the considered data. Despite the small departure from normality with $\mu_1=.05$ and $\mu_2=-.2$, there is a serious impact on the coverage probability of the confidence interval. If dependence is neglected, the interval would not be close to a $95\%$ confidence interval.
\item Case 2: We create a sample of $Y=-\lambda\ln(1-U)$ with $\lambda=0.5, 1, 5, 10, 20, 30$. For each of the values of $\lambda$, we create $N$
samples. Each of the samples produces a confidence interval, and we report proportions of samples that cover the true parameter $\lambda$ as coverage probabilities. For the given value of $\mu_1,$ we have $\sigma^2_f=0.9702667\lambda^2$ and $\sigma_f=0.9953594\lambda$.
\begin{table}[h!]
\centering
\captionsetup{justification=centering,margin=0.5cm}
\begin{tabular}{|c||c|c|c|c |c|c|c|c|c| }
\hline
N/$\lambda$ & 0.5 & 1 & 5 & 10 & 20 & 30\\ 
\hline\hline
100 & 96 &95 &98 &94 &95 &92\\ 
\hline 
100 & 96& 94& 95& 96& 90& 96 \\
\hline
100 &94& 97& 94& 93& 99 & 94 \\
\hline
\end{tabular}
\caption{Coverage probabilities under the Exponential marginal} 
\end{table}
For this example, the low departure from independence doesn't influence much the distribution of the estimator or the confidence interval. Assuming independence in this case would not have serious consequences. When the same sample is used for all values of $\lambda$, we obtain equal coverage probabilities.
\item Case 3. In this case, we check that the Data that shows possible correlation 0, $\rho_S=0$ and $\tau=0$.  We estimate the mean of the marginal distribution and build a confidence interval for samples of 1000 observations. This is repeated 100 times to obtain coverage probabilities. The procedure is repeated for sample sizes 500 and 100.

\begin{table}[h!]
\centering
\captionsetup{justification=centering,margin=0.5cm}
\begin{tabular}{|c||c|c|c|c |c|c|c|c|c| }
\hline\hline
N=1000 & 97 &98 &94 &95 &92 &92\\ 
\hline 
N=500 & 98& 95& 96& 92& 91& 96 \\
\hline
N=100 &95& 95& 91& 96& 94 & 91 \\
\hline
\end{tabular}
\caption{Coverage probabilities under uniform marginal} 
\end{table}

\begin{table}[h!]
\centering
\captionsetup{justification=centering,margin=0.5cm}
\begin{tabular}{|c||c|c|c|c |c|c|c|c|c| }
\hline\hline
N=1000 & 99 &91 &95 &97 &99 &95\\ 
\hline 
N=500 & 95& 94& 97& 97& 91& 97 \\
\hline
N=100 &97& 97& 97& 94& 95 & 98 \\
\hline
\end{tabular}
\caption{Coverage probabilities under uniform marginal and i.i.d} 
\end{table}
These tables also show that the confidence intervals under independence tend to capture more than $94\%$ of the time the true value of the mean. 
\item Case 4: Estimation of $\mu_1$ for a sample with the uniform marginal distribution. We construct  100 confidence intervals for each of the true values of $\mu_1=.05, .1, .11$. The coverage probabilities are reported for sample sizes $N=100, 500, 1000$. These intervals use $w=0.25, 0.5, 0.75$.

\begin{table}[ht!]
\centering
\captionsetup{justification=centering,margin=0.5cm}
\begin{tabular}{|c||c|c|c|c |c|c|c|c|c| }
\hline
$w$/$\mu_1$ & .05 & .05 &.1 &.1 &.1 &.11 & .11 & .11\\
\hline\hline
.25 & 100 &100 &100 &100 &100 &100 & 100& 100\\ 
\hline 
0.5 & 98& 100& 100& 100& 100& 100 &99 & 99\\
\hline
0.9 &95& 93& 95& 95& 94 & 96 & 95& 93\\
\hline
1 & 94 &93 &95 &95 &94 &96 & 95& 93 \\
\hline
\end{tabular}
\caption{Coverage probabilities for $\mu_1$ as function of $w$ and $\mu_1$.} 
\end{table}
This table shows that the best estimator would use values of $w$ close to 1, giving more weigh tho the impact of $\varphi_1(x)$ on the estimator of $\mu_1$. This is true because when $w\le 0.6$, even for a sample of size $5000$, confidence intervals for $\mu_1$ are too wide and cover  the true value of $\mu_1$ almost always.
\end{enumerate}

\section{Conclusion and comments}
We have provided a characterization of copulas with square integrable symmetric densities. We have used this representation to construct several families of copulas. One of the constructed copula families is an extension of the FGM copula family. Other families include trigonometric functions and are all new in the literature. Mixing properties of the constructed copula families have been established. Theorems have been provided for $\rho$-mixing and $\psi$-mixing of Markov chains generated by these copulas and the uniform distribution. This result is equivalent to the said mixing for Markov chains generated by means of a continuous marginal distribution. 

Copulas based on $\{1, \sin2\pi kx, \cos 2\pi kx, k\in\mathbb{N} \}$ have been proposed. Their Spearman $\rho$ and Kendall $\tau$ indicate a way to modify the copula without modifying $\tau(C)=\rho_S(C)=0$. The equalities $\tau(C)=\rho_S(C)=0$ hold for the independence copula. We have investigated central limit theorems for these Markov chains. An example with $\lvert \mu_1\lvert \le .11, \mu_2=-4\mu_1$, $\mu_k=0, k\ge 3$ was used in simulations to illustrate departure from independence, while conserving $\tau(C)=\rho_S(C)=0$. This simulation study has shown that in some cases, the dependence is not seen graphically, and even on the correlation test. However, assuming independence produces poor confidence intervals. Various constructions of this form can be used to model the dependence structure of the data while keeping the key factors of association that the researcher doesn't want to modify. We have provided in this case a central limit theorem for the estimator of $(\mu_1, \mu_2)$; and indicated how it can be used for testing and confidence intervals.

We have shown that examples of copulas based on finite sums of terms can increase $\rho_S(C)$ up to $0.49$ and $\tau(C)$ to $0.32$. These values are larger than those of the popular FGM copula family.  A conclusion on the obtained extension of the FGM copula family is that it can help modify the Kendall coefficient while keeping the Spearman's $\rho$ at the level $\lambda_1$ by introducing non-zero coefficients. This fact is based on $\rho_S(C)=\lambda_1$ for all copulas of the family.

We also mention that functions $\varphi_k(x)$ can act as parameters of the constructed copula families, and be subject to estimation issues. This would be the case when one needs to find the appropriate perturbation that fits the best the relevant situation based on the information at hands.
This question will rely on mixing properties of the constructed copula families and is one of the topics for further research on estimation problems based on these copulas. A drawback for the constructed copula families is that, being absolutely continuous, they do not exhibit any tail dependence. Consideration of their tail dependent extensions is part of ongoing work.

\section{Appendix of proofs}
{\bf Proof of Proposition \ref{pr1}}
\begin{enumerate}
\item Note that for any copula, we have $C(x,1)=x$ and $C(x,0)=0$. Therefore, for any absolutely continuous copula, $\varphi(x)=1$ satisfies
\[  K \varphi(x) = \int_{0}^1c(x,y)dy=C_{,1}(x,1)-C_{,1}(x,0)=1=\varphi(x).\]
It follows that the function $\varphi(x)=1$ is an eigenfunction of any kernel operator defined by an absolutely continuous copula associated to the  eigenvalue $1$. 
\item The proof of the second point relies on formula \eqref{cond1} and the first point of Proposition \ref{pr1}.
\item The proof of this third point is a consequence of the fact that when the eigenvalue $1$ has multipilicity higher than $1$, it generates a subspace of dimension greater than $1$. In this subspace, an orthonormal basis not containing $\varphi(x)=1$ can be constructed. 
$\square$
\end{enumerate}
{\bf Proof of Theorem \ref{T1}}

The statement of Theorem \ref{T1} implicitly claims that the sum has point-wise convergence. To prove this, we recall Weierstrass' M-test.  
By Weierstrass’ M-test, for a sequence $f_n(x)$ defined on the same support $E$, if $\lvert{}f_n(x)\lvert{}\le M_n$ for all $x$ and $\sum M_n<\infty$, then $\sum f_n(x)$ converges uniformly and absolutely on $E$. For the sum \eqref{cop}, Condition \eqref{cond2} is equivalent to 
$$ \sum_{k=1}^{\infty} \lvert{}\lambda_k \alpha_k\lvert{} \le 1.$$
Note that $ \lvert{}\lambda_k \varphi_k(x)\varphi_k(y)\lvert{}\leq \lvert{}\lambda_k\alpha_k\lvert{}$ for all $x,y$ and $k$. Moreover, $\sum_{k=1}^{\infty} \lvert{}\lambda_k \alpha_k\lvert{} $ converges as a consequence of Condition \eqref{cond2}. Thus, by Weierstrass M-test, the series that defines the density of the copula converges absolutely and uniformly. Moreover, \eqref{cop} also holds with absolute convergence and uniform convergence.
The rest of the proof relies on the fact that functions form an orthonormal system. Orthogonality implies $C(1,x)=C(x,1)=x$, because $\int_0^1\varphi_k(s)ds=0$ and $C(x,0)=C(0,x)=0$ by definition. $\square$

{\bf Proof of Theorem \ref{theoM}}

The first condition of this theorem implies $\psi'$-mixing as a consequence of Longla (2015) \cite{Lo2015}. Longla (2015)\cite{Lo2015} showed that if the density of the absolutely continuous part of $C^{s_1}(u,v)$ is bounded away from $0$ on a set of Lebesgue measure 1, then $\psi'_{s_1}>0$. Moreover, $c^{s_2}(u,v)<M$ implies $\psi_{s_2}<\infty$. Therefore, by Bradley (2005) \cite{Bradley}, the first condition of Theorem \ref{theoM} implies $\psi$-mixing.
The second condition implies $\psi$ as a result of Longla et al. (2022c) \cite{LoMous}. $\square$

{\bf Proof of Example \ref{rem1}. }

If this copula generates $\psi$-mixing stationary Markov chains with continuous marginal distributions, then it would have to generate $\psi'$-mixing and would therefore generate $\rho$-mixing (contradiction) with Longla (2014)\cite{L2}. Therefore, it doesn't generate $\psi$-mixing. Though its density is $c(u,v)=0$ for all $(u,v)\in (0,1)^2$ such that $v\ne u\ne 1-v$, this density doesn't exist on a set of Lebesgue measure $0$ and is bounded on a set of Lebesgue measure $1$. The second example is based on the fact that the density of the given copula is equal to $2$ on a set of Lebesgue measure $1/2$ and equal to $0$ on its complement. $\square$

{\bf Proof of Lemma \ref{sincosl}}

The upper bound on these values is based on the assumptions on the coefficients of the series. These conditions imply that $\lvert{}\mu_k\lvert{}\le 1/2$ and the largest value in absolute value is achieved when $\lvert{}\mu_1\lvert{}=1/2$, $\mu_k=0$ for $k>1$ and $\lambda_k=0$ for $k\ge 1$. Moreover, \[ \lvert{}\sum_{k=1}^{\infty} \frac{\mu_k}{k^2}\lvert{}\leq \sum_{k=1}^{\infty} \frac{\lvert{}\mu_k\lvert{}}{k^2}\leq \sum_{k=1}^{\infty} \lvert{}\mu_k\lvert{}\leq 1/2 \quad \text{and}\] 
\[ \lvert{}\sum_{k=1}^{\infty}\frac{4\mu_k+2\lambda_k\mu_k}{k^2}\lvert{}\leq 2(2+\lambda_1)\lvert{}\mu_1\lvert{}+\sum_{k=2}^{\infty}2\frac{(2+\lambda_k)\lvert{}\mu_k\lvert{}}{k^2} \]\[\le 2(2+\lambda_1)\lvert{}\mu_1\lvert{}+\sum_{k=2}^{\infty}2\lvert{}\mu_k\lvert{}=2(1+\lambda_1)\lvert{}\mu_1\lvert{}+2\sum_{k=1}^{\infty}\lvert{}\mu_k\lvert{}. \]

\[\text{Therefore,}\quad \lvert{}\sum_{k=1}^{\infty}\frac{4\mu_k+2\lambda_k\mu_k}{k^2}\lvert{}\leq2(1+\lambda_1)\lvert{}\mu_1\lvert{}+1\leq 2.\]
The last inequality uses $\lvert{}\lambda_1\lvert{}+\lvert{}\mu_1\lvert{}\leq 1/2$. Equality happens for $\lvert{}\mu_1\lvert{}=1/2$ and $\lambda_1=0$.  These inequalities justify the boundary on $\rho_S(C)$ and $\tau(C)$.
$\square$

{\bf Proof of Theorem \ref{mixcop}}

The proof of Theorem \ref{mixcop} relies on the fact that $C^n(u,v)$ is the $n^{th}$-fold product of $C(u,v)$ , $\varphi_k(x)$ form an orthonormal basis of $L^2(0,1)$, $\lambda^n_k$ are all eigenvalues of the Hilbert-Schmidt operator associated to $c^n(u,v)$.
Note that though we have square integrability, the supremum in Theorem \ref{mixcop} can be equal to $1$. The rest is an application of Longla and Peligrad (2012) \cite{LoMag}. $\square$.

{\bf Proof of Theorem \ref{decompPsi}}

The proof of Theorem \ref{decompPsi} is an application of Longla et al (2022c) \cite{LoMous}. In Longla et al (2022c) \cite{LoMous} it was shown that if the density $c^n(u,v)< 2$ on $[0,1]$ for some $n$, then the copula $C(u,v)$ generates $\psi^*$-mixing. Let $\underline{\lambda}=\sup_{k}\lvert{}\lambda_k\lvert{}<1$. From the fact that the density is square integrable, we have \[ c^n(x,y)\ge 1-\sum\lvert{}\lambda_k\lvert{}^n(2k+1)\alpha_k>\] \[1-\underline{\lambda}^{n-3}\sum\lvert{}\lambda_k\lvert{}^3 (2k+1)\alpha_k>1-M\underline{\lambda}^{n-3}\ge 0.\]  Therefore, $c^n(u,v)$ is bounded away from 0 on a set of Lebesgue measure $1$. This implies $\psi'$-mixing.  Moreover, $\lvert{}\lambda_{k}\lvert{}\le (2k+1)^{-1}$ when $k$ is even. Thus, \[ c^n(x,y)\le 1+\sum\lvert{}\lambda_k\lvert{}^n(2k+1)\le 1+\sum_{odd}\lvert{}\lambda\lvert{}_k^{n-1}+\sum_{even}\lvert{}\lambda_k\lvert{}^n (2k+1)\le\] \[ 1+\underline{\lambda}^{n-3}\sum_{odd}\lambda_k^2+\sum_{even}\lvert{}\lambda_k\lvert{}^n (2k+1).\]
The last inequality uses the fact that the sequence $\lambda_k$ converges to $0$ or has finitely many values. The second portion is also bounded. It is written separately to emphasize the fact that for even values of $k$, it is possible to have larger values when $\lvert{}\lambda_k\lvert{}>(2k+1)^{-1}$. For odd values of $k$, this is not possible because $\lvert{}\lambda_k\lvert{}\leq (2k+1)^{-1}$. If we denote $M=\sum_{k=1}^\infty\lambda_k^2$, then $M<\infty$ because $c(x,y)$ is square integrable. So, we obtain
\[1-M\underline{\lambda}^{n-3}\le c^{n}(x,y)\leq 1+M\underline{\lambda}^{n-3}.\]
$M$ is a constant free of $n$. So, we conclude that as $n\to\infty$, $M\underline{\lambda}^{n-3}\to 0$. Therefore, $c^n(x,y)<2$ for sufficiently large values of $n$. By Longla et al (2022c) \cite{LoMous} it follows that the Markov chains are $\psi$-mixing. 

Note that there are still some values of $\lambda_k>(2k+1)^{-1}$ that are not considered by Theorem \ref{decompPsi} which can increase up to the minimum of the considered even Legendre polynomials. Using the same arguments, we extend the theorem to the cases when the supremum of $\lambda_k$ is less than $1$. Due to square integrability, we have $\sum\lambda_k^2<\infty$. The comparison test for series leads to $\lambda_k^2(2k+1)\to 0$ as $k\to \infty$ because $\sum(2k+1)^{-1}$ diverges. Therefore, there exists an integer K, such that for $k>K$, $\lambda^2_k\leq \varepsilon/(2k+1)$. Now, for $k\leq K$, it is easy to find $N$ such that for all $n\ge N$, $\sum \lambda_k^n(2k+1)<\varepsilon_1$, where the sum is taken over even integers less than or equal to $K$.
This leads to \[ c^n(u,v)\leq 1+\underline{\lambda}^{n-3}\sum_{odd}\lambda^2_k+\] \[\underline{\lambda}^{n-N}\sum_{even k\leq K}\lvert{}\lambda_k\lvert{}^N(2k+1)+\underline{\lambda}^{n-4}\sum_{even k>K}\lvert{}\lambda_k\lvert{}^4(2k+1).\]
\[\text{So,}\quad c^{n}(u,v)\leq 1+\underline{\lambda}^{n-3}M+\underline{\lambda}^{n-N}\varepsilon_1+\underline{\lambda}^{n-3}M\varepsilon.\]
The last inequality has parameters $\varepsilon, \varepsilon_1, N$ and $M$ that are free of $n$. Therefore, as a sequence that converges to 1, the right hand side can be made strictly less than 2 for any values of these parameters as $n$ gets sufficiently large. So, $c^{n}(u,v)<2$ for some integer $n$. 
$\square$

{\bf Proofs of Theorem \ref{mut} and Corollary \ref{colmu}}

Assume $(U_1,\cdots, U_n)$ is a Markov chain with uniform marginals generated by copula \eqref{model}. 

{\bf Proof of Theorem \ref{mut}}

Let $f_i(u,v)=\varphi_i(u)\varphi_i(v)$, for $i=1,2$ and $f=af_1+bf_2$. It is obvious that $\mathbb{E}f(U_i,U_{i+1})=a\mu_1+b\mu_2$ and  $var(f(U_1,U_2))$ is finite. Moreover, $Y_{i-1}=(U_{i-1},U_{i}),$ $i=2\cdots n$ is a Markov chain. This Markov chain is ergodic because the original Markov chain is $\psi$-mixing. 
Therefore, by Kipnis and Varadhan (1986)\cite{KV}, the CLT holds with 
$$\sigma^2_f=var(f(U_1,U_2))+2\sum_{i=1}^{\infty}cov(f(U_1,U_2), f(U_{1+i},U_{2+i})).$$

Note that for a Markov chain generated by \eqref{model},  the density of the cumulative distribution function of $(U_1,U_2, U_{1+i}, U_{2+i})$ is obtained using the Markov property as follows. For $i\ge 2$, using formula \eqref{joint}, $(U_2, U_{1+i})$ has density 
$$h_1 (u,v)=1+\mu^{i-1}_1 \varphi_1(u)\varphi_1(v)+\mu^{i-1}_2 \varphi_2(u)\varphi_2(v).$$
Therefore, the joint density of $(U_1,U_2, U_{1+i}, U_{2+i})$ is 

\begin{equation}
h(s, u,v, w)=c(s,u)h_1(u,v)c(v,w). \label{cf}
\end{equation}
 Formula \eqref{cf} implies 
$$\mathbb{E}(f_1^2(U_1,U_2))=1+\mu_1\large\left(\int_0^1\varphi_1^3(t)dt\large\right)^2+\mu_2\large\left(\int_0^1\varphi_2(t)\varphi_1^2(t)dt\large\right)^2, $$
$$\mathbb{E}(f_1(U_1,U_2) f_1(U_{2},U_{3}))=\mu_1^2\int_0^1\varphi_1^4(t)dt,$$
$$\mathbb{E}(f_1(U_1,U_2)f_1(U_{1+i},U_{2+i}))=\mu_1^2+\mu^{1+i}_1\large\left(\int_0^1\varphi_1^3(t)dt\large\right)^2, \quad i>1;$$
$$\mathbb{E}(f_1(U_1,U_2)f_2(U_1,U_2))=\mu_1\large\left(\int_0^1\varphi_1^2(t)\varphi_2(t)dt\large\right)^2+\mu_2\large\left(\int_0^1\varphi_1(t)\varphi_2^2(t)dt\large\right)^2, $$
$$\mathbb{E}(f_1(U_1,U_2) f_2(U_{2},U_{3}))=\mu_1\mu_2\int_0^1\varphi_1^2(t)\varphi_2^2(t)dt,$$
$$\mathbb{E}(f_1(U_1,U_2) f_2(U_{1+i},U_{2+i}))=\mu_1\mu_2^i\int_0^1\varphi_1^2(t)\varphi_2(t)dt\int_0^1\varphi_2^3(t)dt+$$
$$+\mu_2\mu_1^i\int_0^1\varphi_2^2(t)\varphi_1(t)dt\int_0^1\varphi_1^3(t)dt+\mu_1\mu_2.$$
Corresponding formulas for $f_2$ are not provided here. They are obtained by symmetry. Using  $\varphi_1(x)=\sqrt{2}\sin(2\pi x)$ and $\varphi_2(x)=\sqrt{2}\sin(4\pi x)$, we obtain $$\int_0^1\varphi^3_i(x)dx=0, \quad \int_0^1\varphi^2_1(x)\varphi^2_2(x)dx=1$$ $$\text{and}\quad \int_0^1\varphi_1^2(x)\varphi_2(x)dx=\int_0^1\varphi_1(x)\varphi^2_2(x)dx=0, \quad \int_0^1\varphi_i^4(x)dx=\frac{3}{2}.$$ Therefore, substituting and simplifying leads to
$\sigma^2_f=a^2+b^2-2ab\mu_1\mu_2.$
It follows that for $\bar{f_n}=\frac{1}{n-1}\sum_{i=1}^{n-1}(a\varphi_1(U_i)\varphi_1(U_{i+1})+b\varphi_2(U_i)\varphi_2(U_{i+1}))$,  \begin{equation}\sqrt{n-1}(\bar{f_n}-a\mu_1-b\mu_2) \to N(0, a^2+b^2-2ab\mu_1\mu_2). \label{last}\end{equation}
Therefore, the Cramer-Wold device concludes the proof of Theorem \ref{mut}.  
$\square$

{\bf Proof of Corollary \ref{colmu}}

 According to our work, this Markov chain is $\psi$-mixing.  Let $f(x)$ be a function such that $\mu=\int_0^1f(x)dx$ and $\int_0^1f^2(x)dx<\infty$. Define $S_n(f)=\sum_{i=1}^n f(U_i)$. By Kipnis and Varadhan (1986)\cite{KV}, the central limit theorem holds in the form 
\begin{equation}
\sqrt{n}(\frac{S_n(f)}{n}-\mu) \to N(0, \sigma^2_f), \quad \text{where},\quad \sigma^2_f=\lim_{n\to\infty} \frac{var(S_n(f))}{n}.
\end{equation}
Based on the provided formula of joint cumulative distribution of $(U_0, U_n)$ and stationarity,  $\varphi_1(x)=\sqrt{2}\sin 2\pi x$ and $\varphi_2(x)=\sqrt{2}\sin 4\pi x$ imply $$cov(f(U_1), f(U_k))=\mu_1^{k-1}\int_0^1\int_0^1 \varphi_1(u)\varphi_1(v)f(u)f(v)dudv+$$ $$\mu_2^{k-1}\int_0^1\int_0^1 \varphi_2(u)\varphi_2(v)f(u)f(v)dudv.$$ 
Moreover, in the context of reversible Markov chains, we have \begin{equation*}\sigma^2_f=\sigma^2+2\sum_{k=2}^{\infty}cov(f(U_1), f(U_k))=\sigma^2+2(\frac{\mu_1A_1^2}{1-\mu_1}+\frac{\mu_2A_2^2}{1-\mu_2}),\end{equation*}
where $A_i=\int_0^1\varphi_i(u)f(u)du$ and $\sigma^2$ is the variance of $f(X_i)$. In this case, $\mu_2=-4\mu_1$. Simple computations lead to the formula of $\sigma^2_f$.
\begin{itemize}
\item Case 1: Dependent Bernoulli observations. Select $f(x)=\mathbb{I}(x\le a)$ for $a\in(0,1)$. Via simple computations, we obtain $\mu=\mathbb{E}f(U_1)=a$, $\sigma^2=a(1-a)$ and $A_i^2=\frac{(1-\cos 2 \pi i a)^2}{2\pi^2i^2}=\frac{2\sin^4 \pi i a}{\pi^2i^2}$. Thus, equation \eqref{I} holds.

\item Case 2: Trigonometically dependent exponential sequences.  Consider $f(x)=-\lambda\ln(1-x)$. This defines a sequence of variables $X_1,\cdots, X_n$ that has exponential marginal distributions ($X_i=-\lambda\ln (1-U_i)$). $\mu=\mathbb{E}(X_i)=\lambda$, $\sigma^2=\lambda^2$. Simple computations give  formula \eqref{I2}

\item Case 3. Dependent uniform Data. Consider now $f(x)=x$. The generated Markov chain itself. $\mathbb{E}(X)=1/2$, $\sigma^2=1/12$. $A_i^2=\frac{1}{2i^2\pi^2}$. So, formula \eqref{vm} holds. For $\lvert \mu_1\lvert \le .11$, the variance \eqref{vm} is strictly positive. 
\item Estimating $\mu_1$. Consider estimating the parameter $\mu_1$ based on the generated Markov chain. Take any $w\in\mathbb{R}$ and
$$ f(x,y)=2w \sin(2\pi x)\sin(2\pi y)-\frac{1}{2}(1-w)\sin(4\pi x)\sin(4\pi y).$$
Theorem \ref{mut} and Formula \eqref{last} conclude the proof.
$\square$
\end{itemize}

\section*{Data availability and conflict of interest}
This manuscript has no associated data and there is no conflict of interest.

\bibliographystyle{amsplain}

\end{document}